\documentclass[12pt]{article}
\usepackage{hyperref}
\newcommand{\nocopyright}{
No Copyright\thanks{
The author(s) hereby waive all copyright
and related or neighboring rights to this work,
and dedicate it to the public domain.
This applies worldwide.
}}
\title{Conway's drum quilts}
\author{
Peter G. Doyle}
\date{Version 1.0 dated 15 June 2020
\\ \nocopyright
}
\usepackage{amsthm,amssymb}
\usepackage{verbatim}
\usepackage{amsmath}
\usepackage{graphicx}

\newcommand{\putfigwithsize}[4]{
\begin{figure}[btp]
\centerline{\mbox{\includegraphics[#1]{figures/#3.pdf}}}
\caption{#4}
\label{fig:#2}
\end{figure}
}

\newcommand{\figref}[1]{\ref{fig:#1}}

\newcommand{\cross}{\times}

\newcommand{\fig}[1]{}

\begin{document}
\maketitle

\begin{abstract}
A `transplantable pair' is a pair of glueing diagrams that can be used to
create pairs of plane domains 
that are isospectral for the Laplace operator.
We present a host of transplantable pairs
worked out by John Conway
using his theory of quilts.
\end{abstract}

\section{Introduction}
A `transplantable pair' is a pair of glueing diagrams that can be used to
create pairs of plane domains or other spaces that
are isospectral for the Laplace operator.
The pair is metonymically `transplantable' because isospectrality of the
glued spaces can be proven using Peter Buser's transplantation
method,
as explained by Buser, Conway, et al.\ \cite{bcds:drum} and Conway \cite{conway:sensual}.

John Conway produced a host of transplantable pairs by applying his
theory of quilts to small projective groups.
I helped by watching admiringly.
In this paper I have reproduced his catalog of pairs.
The sizes of these pairs are 7, 11, 13, 15, and 21.
In \cite{bcds:drum} we presented
the sixteen pairs of sizes 7, 13, and 15,
which are treelike and thus give planar isospectral domains.
We also presented one of the size 21 pairs,
here labeled pair $21(7)$, which yields the `homophonic' domains
shown in Figure \figref{peacocks}.
Conway 
\cite[p. 249]{conway:things}
dubbed these domains `peacocks rampant and couchant'.

\putfigwithsize{width=12cm}{peacocks}{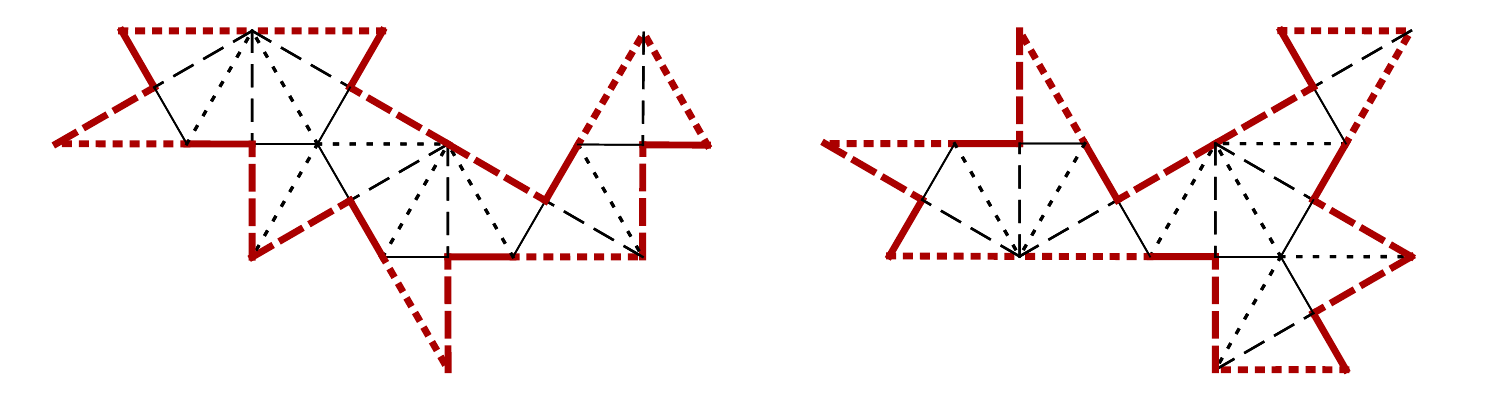}{Peacocks rampant and couchant}

A transplantable pair can be thought of as arising from a pair
of finite permutation actions of the free group on generators $a,b,c$.
These actions are equivalent as linear representations but (if
the pair is to be of any use) not as permutation representations.
In the examples at hand, each representing permutation is an involution,
so these representations factor through the
quotient $F=\langle a,b,c:a^2=b^2=c^2=1 \rangle$,
the free product of three copies of the group of order $2$.

From any transplantable pair
we can get other pairs through the process of
braiding, which amounts to precomposing the permutation representations
with automorphisms of $F$,
called $L$ and $R$:
\[
L:(a,b,c) \mapsto (aba^{-1},a,c)
;
\]
\[
R:(a,b,c) \mapsto (a,c,cbc^{-1})
.
\]
Left-braiding a permutation representation $\rho$
can be viewed as first conjugating $\rho(b)$ by $\rho(a)$, i.e.\ applying
the permutation $\rho(a)$ to each index in the cycle representation of
$\rho(b)$, and then switching $\rho(a)$ with the new $\rho(b)$.
Right-braiding is the same, only with $c$ taking over the role of $a$.

Pairs of permutations that are equivalent in this way belong to the same \emph{quilt}.
We identify pairs that differ only by permuting $a,b,c$, or by 
reversing the pair.
A quilt has extra structure which we are ignoring here:
See Conway and Hsu \cite{conwayhsu:quilts}.
This structure makes it easier to understand and enumerate the pairs.
But the computer has no trouble churning out all the
pairs belonging to the same
quilt.

So despite the title and what you might reasonably expect,
the only place you will find quilts here is in Appendix \ref{appendix},
which reproduces Conway's original quilt calculations.

First we will present the glueing diagrams, and then the
corresponding hyperbolic orbifolds.

\section{Transplantation diagrams}

In the diagrams that follow,
the points being permuted are represented by triangles.
The permutations corresponding to $a,b,c$ are
represented by lines of three styles:
dotted, dashed, and solid.
Black lines separate pairs of points that are interchanged
by the permutation, while
red lines (which are also made thicker) indicate fixed points.
Red lines in the interior of the diagram separate points that
are each fixed, rather than interchanged.
Sometimes black lines occur on the boundary of the diagram,
which means that the boundary must be glued up.
The computer has taken care to lay out the diagrams so that there is at most
one pair of thin boundary lines of each type (dotted, dashed, or solid),
so that even without the usual glueing arrows
there is no ambiguity of how the boundary is to be glued up.

To refer to these pairs, we will write $7(1)$ for the first pair of quilt 7,
$13a(5)$ for the fifth pair of quilt $13a$, etc.
This numbering is canonical, given the starting triple $(a,b,c)$, because
the quilt has been explored by a `left-first search'.
More canonical, but more cumbersome, are the Conway symbols of the
simplest hyperbolic orbifolds that can be obtained from the pair:
See Section \ref{hype}.

For historical reasons, the four quilts of size 11 are called
$11f,11g,11h,11i$.
The missing quilts (whose diagrams actually have size 12)
are to be found in raw form in Appendix \ref{appendix}.

\clearpage
\newcommand{\quilt}[1]{\putfigwithsize{width=\textwidth}{#1}{fig#1}{Quilt #1}}
\newcommand{\quiltsqueeze}[1]{\putfigwithsize{width=\textwidth}{#1}{fig#1}{Quilt #1}}
\quilt{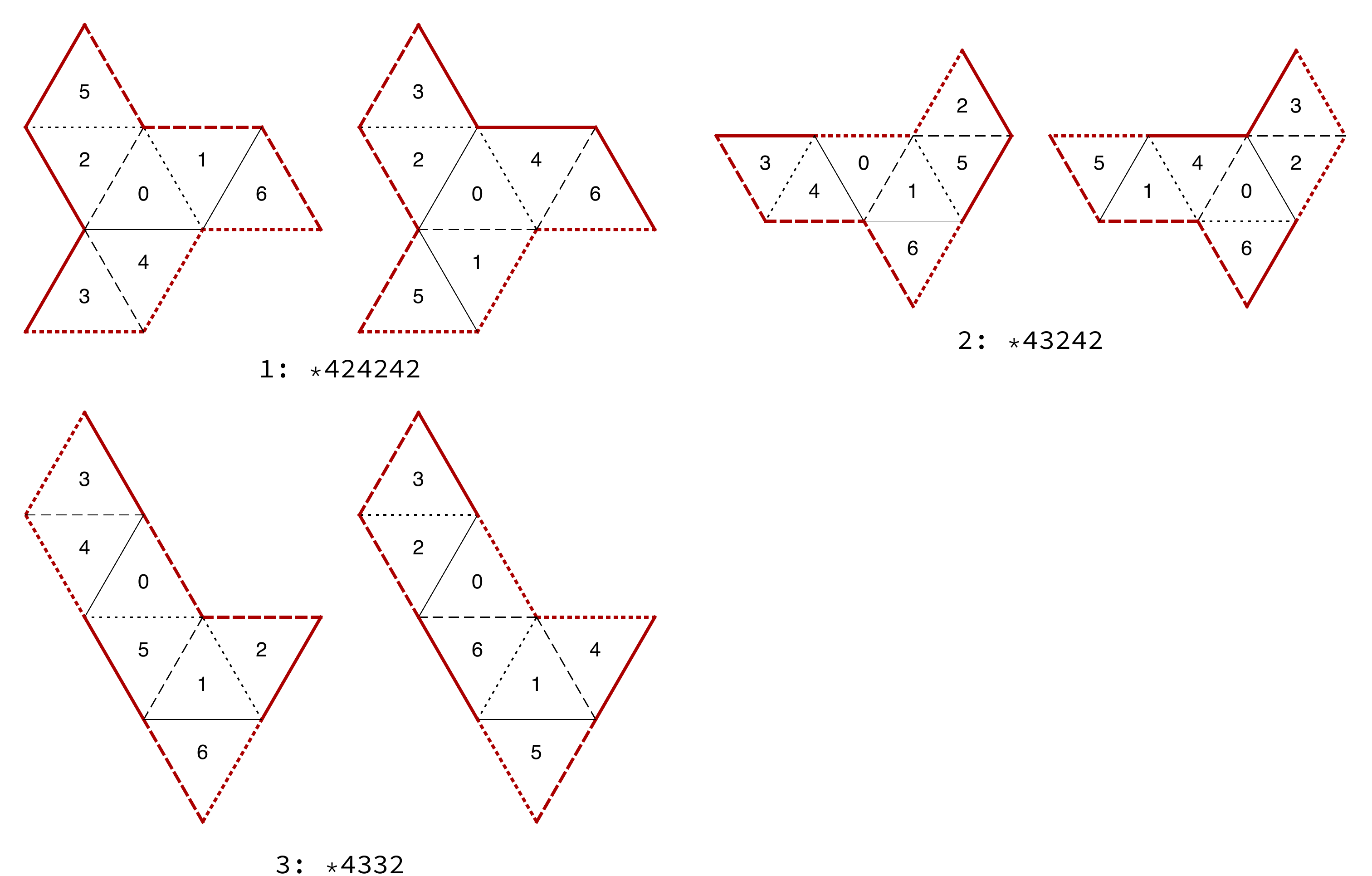}
\quilt{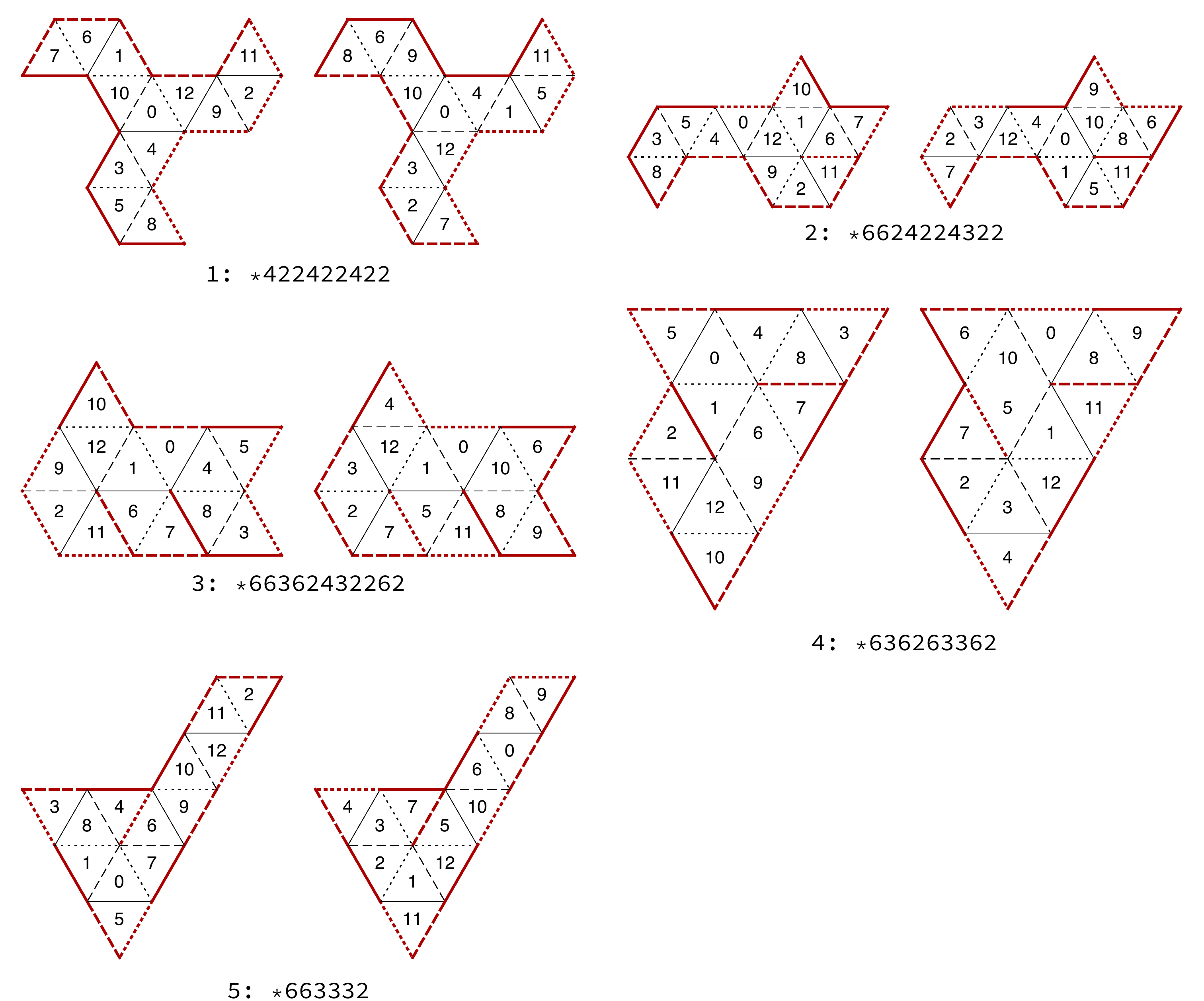}
\quilt{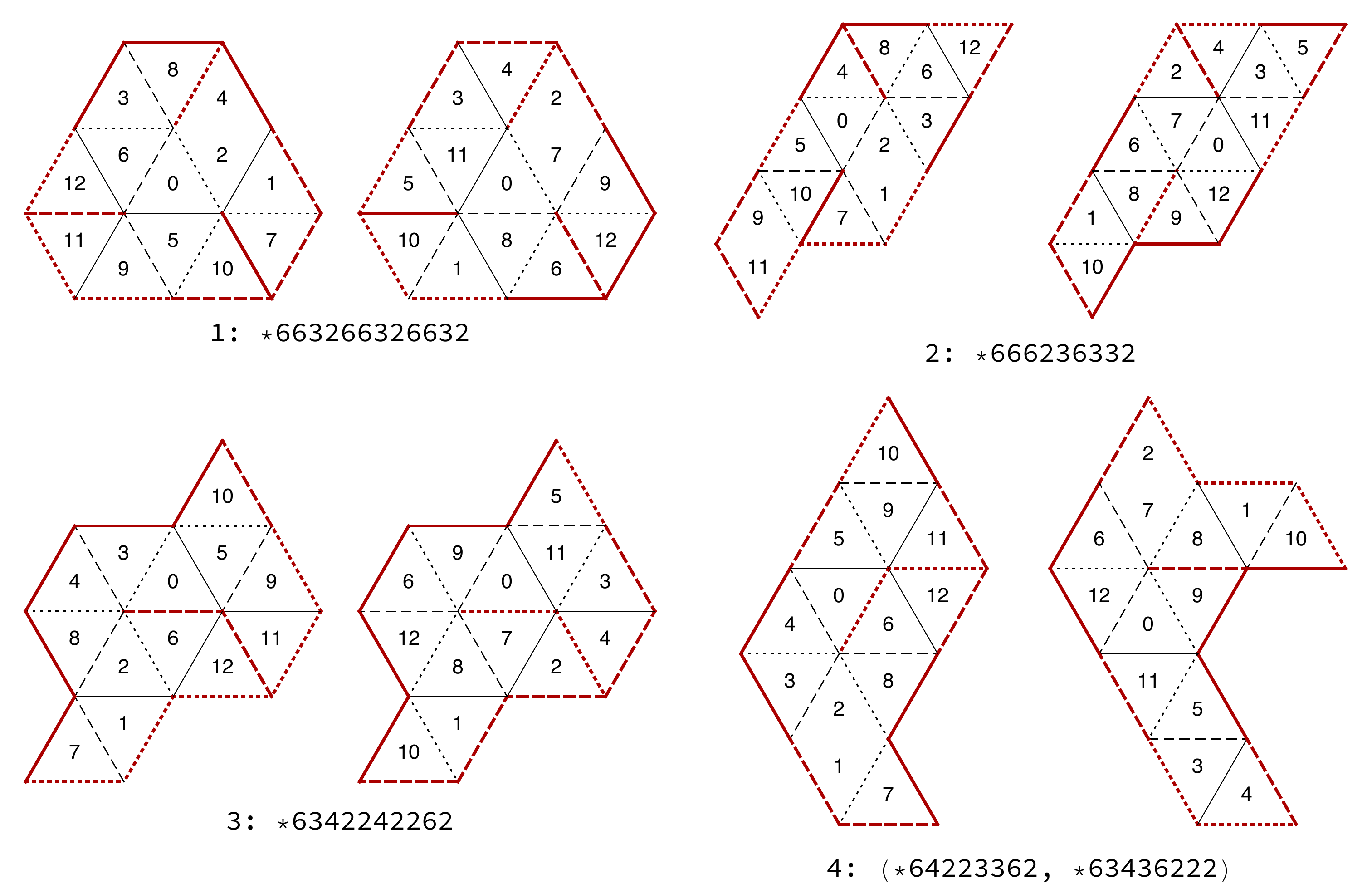}
\quilt{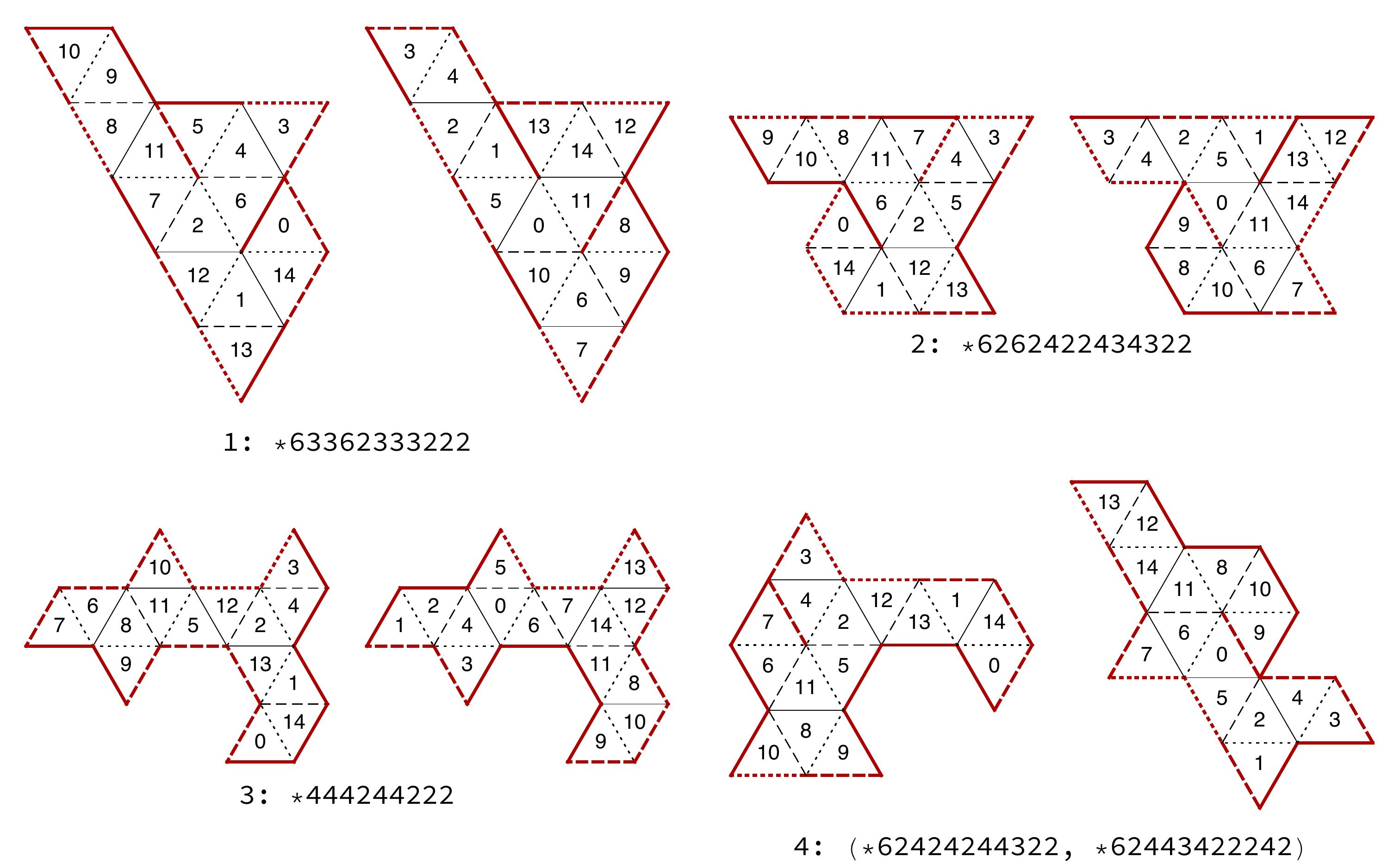}
\quiltsqueeze{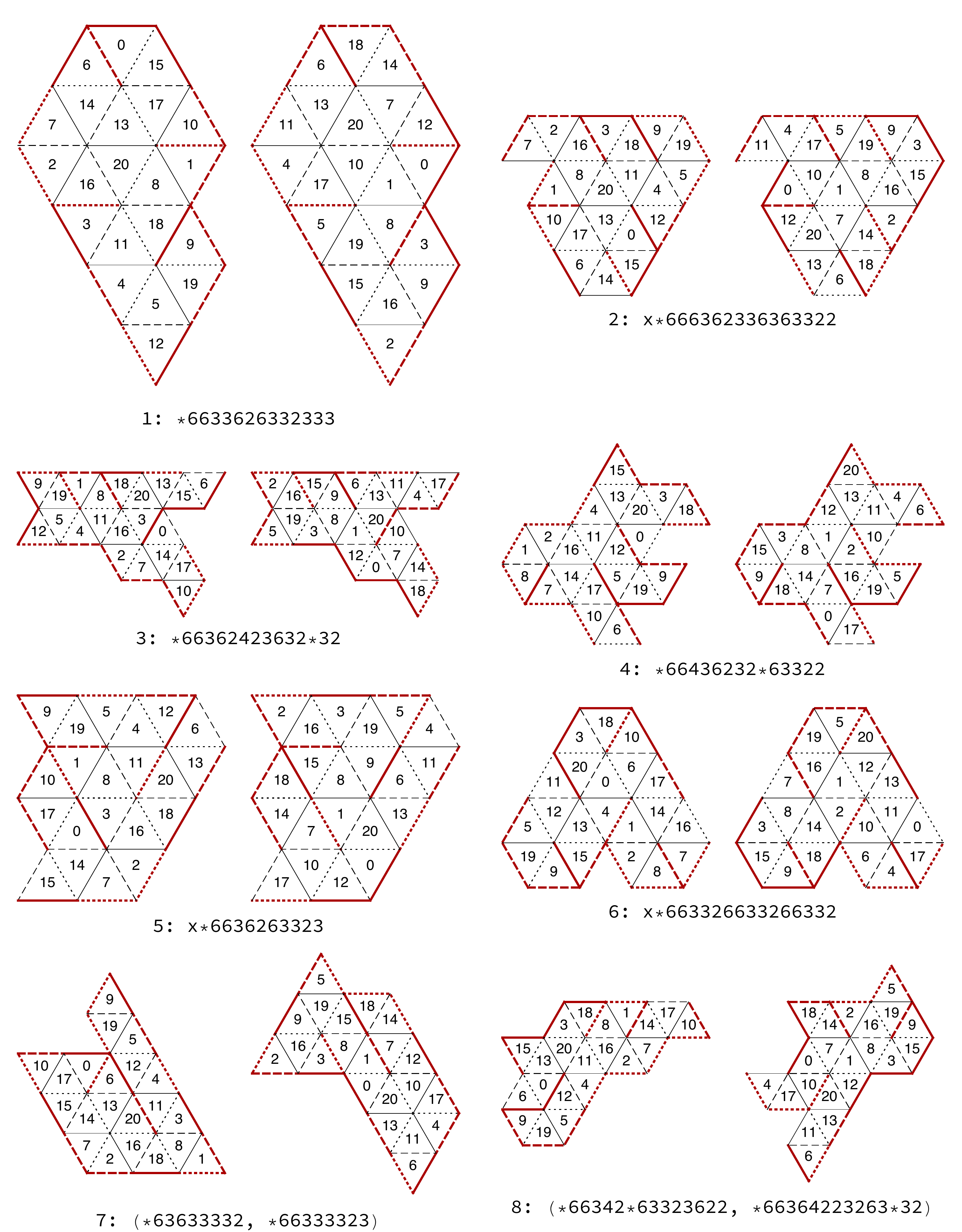}
\clearpage
\quilt{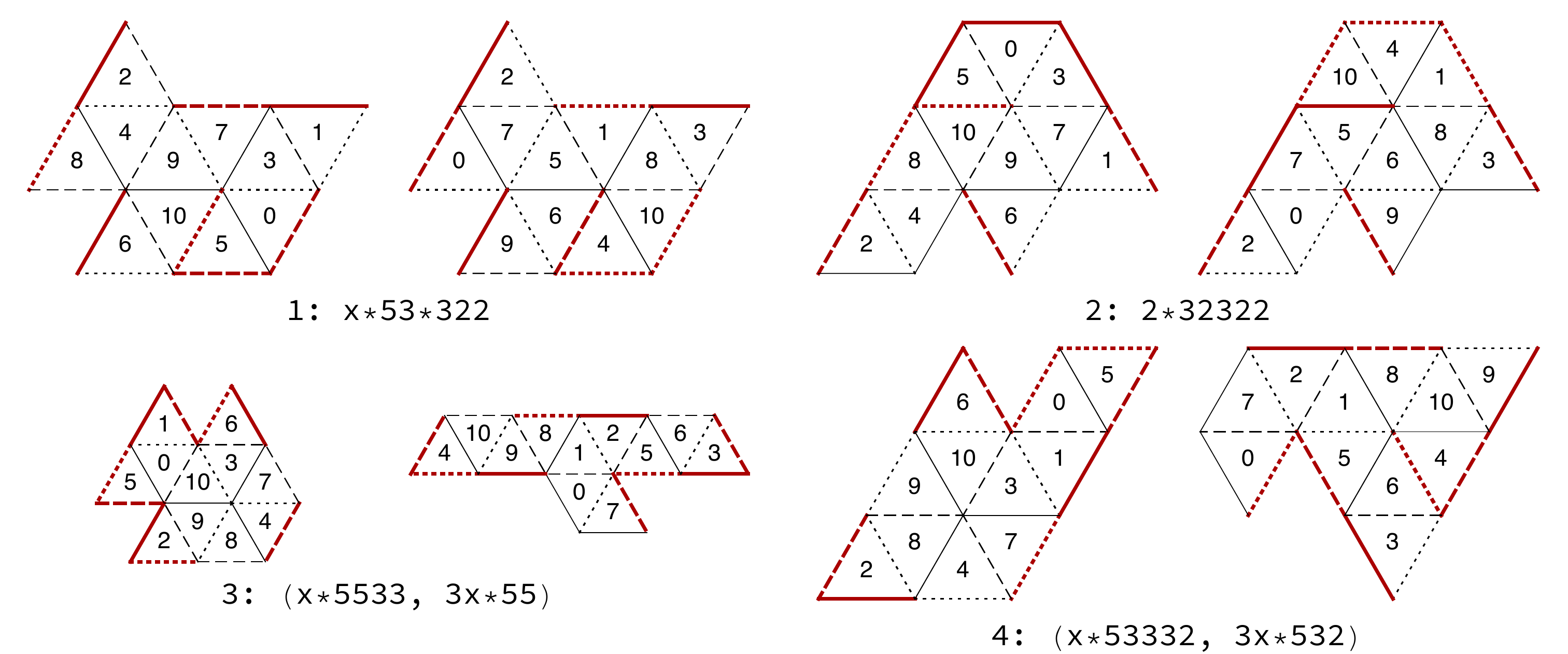}
\quilt{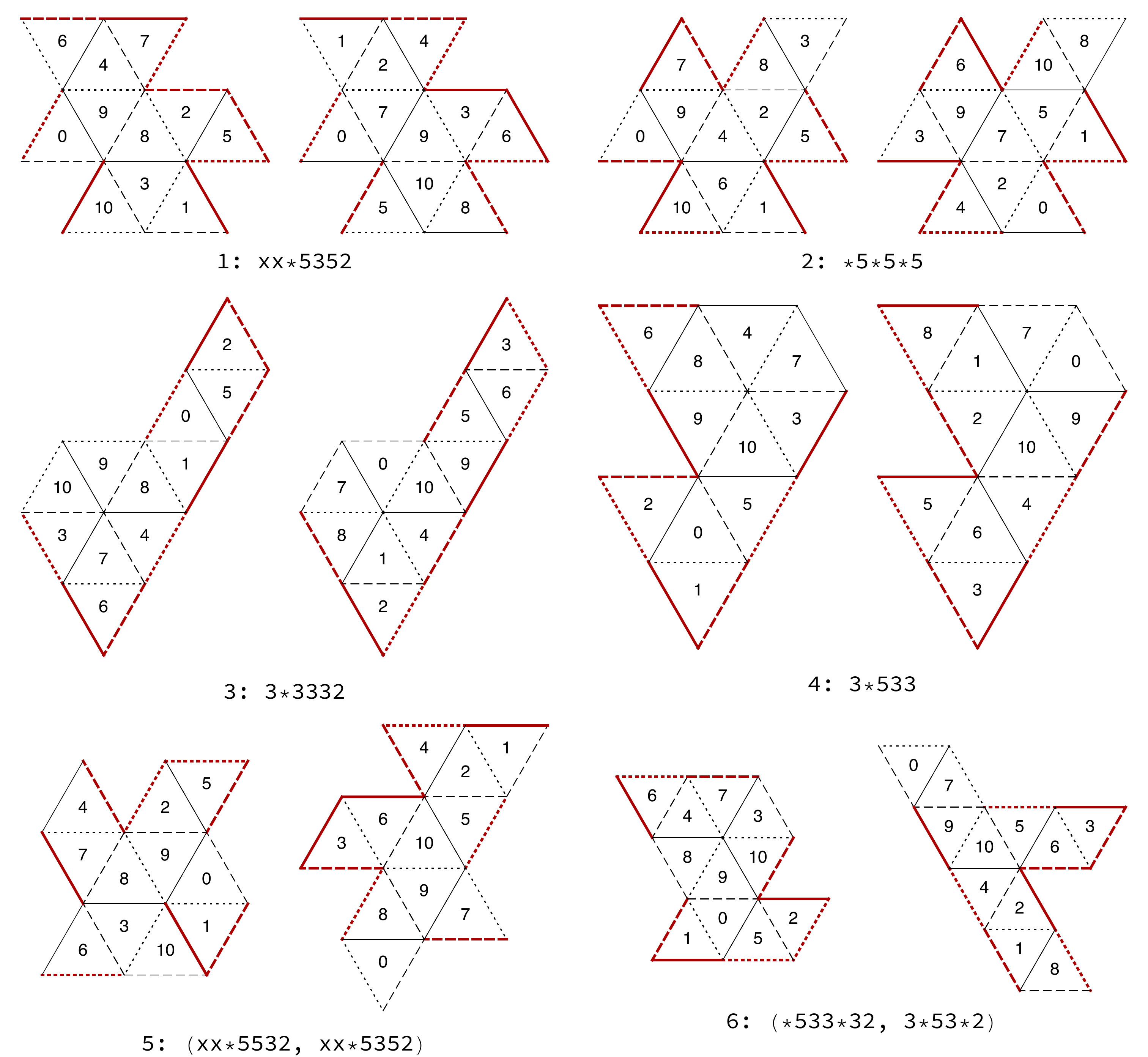}
\quilt{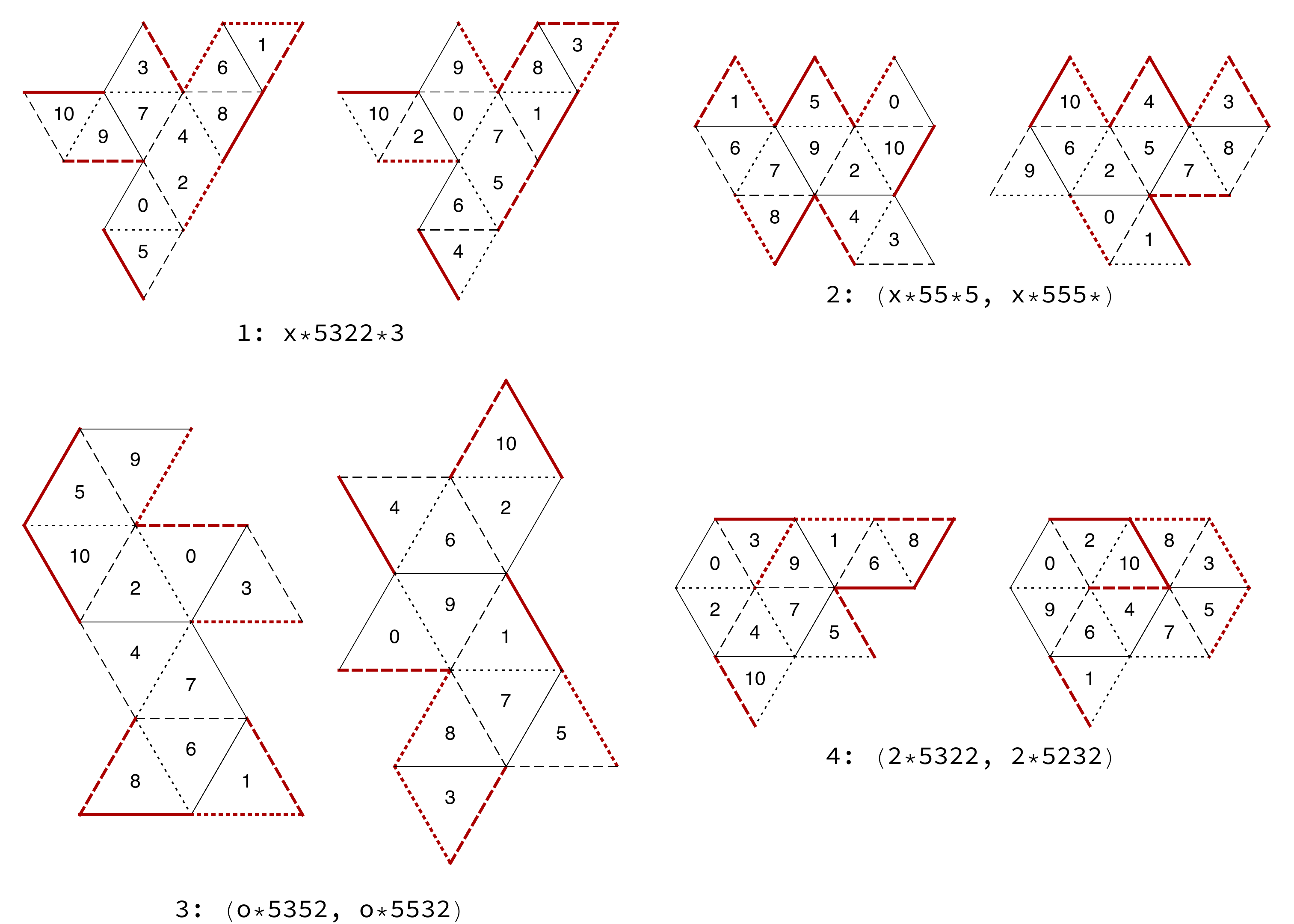}
\quilt{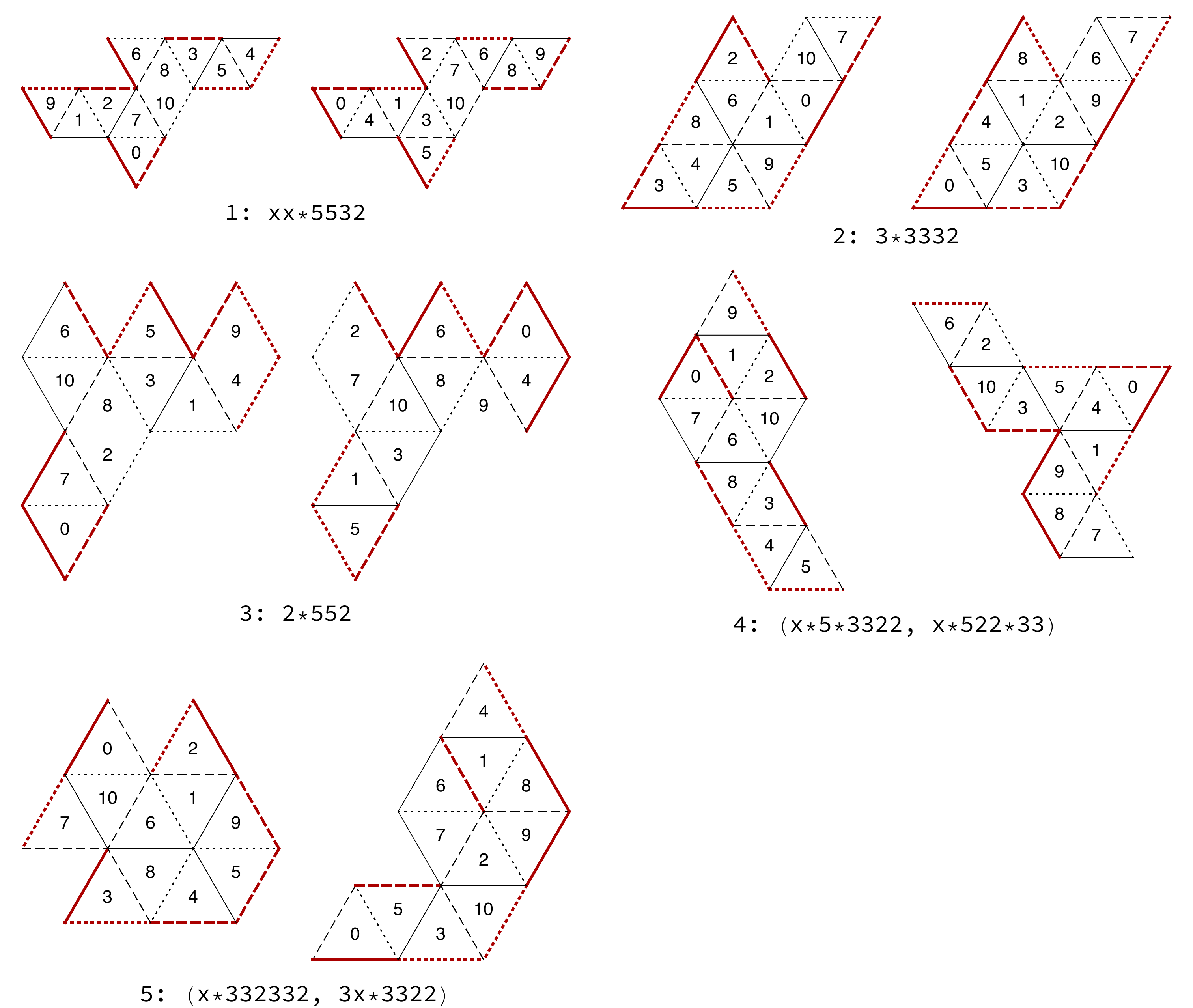}

\clearpage
\section{Hyperbolic orbifolds} \label{hype}

\putfigwithsize{width=.75\textwidth}{sunada6}{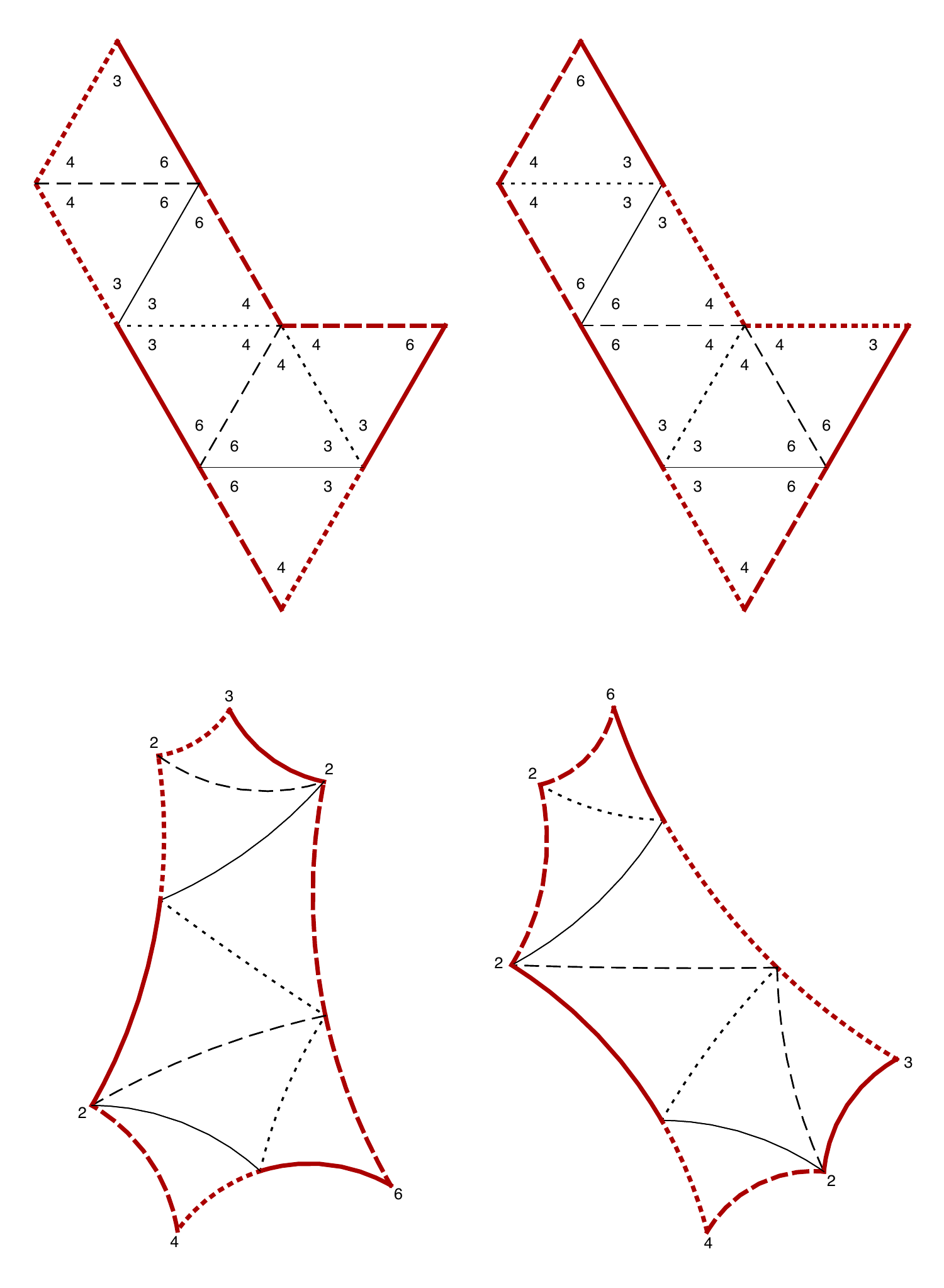}{Isospectral
hyperbolic hexagons arising from diagram $7(3)$.}

Here are the hyperbolic orbifolds corresponding to these transplantable
pairs.
These are the simplest hyperbolic orbifolds that can be produced using the
glueing data.
To get them, for the basic triangle we prescribe angles
just small enough to make each interior cone point
have cone angle evenly dividing $\tau$,
and each boundary corner have angle evenly dividing $\tau/2$.
When the members of the pair differ only by permuting the labels $a,b,c$,
the resulting orbifolds are isometric:
To get non-isometric pairs we will need to destroy this symmetry
by taking one or more of the triangles smaller.
Thus, for example, from the pair $7(3)$ we get the isospectral hexagons shown
in Figure \figref{sunada6}.
This is presumably the simplest pair of isospectral hyperbolic 2-orbifolds.

Each figure gives Conway's notation
(see Conway \cite{conway:notation})
for the associated
orbifolds.

{\bf Note.}
The pairs show here are isospectral as hyperbolic
2-orbifolds.
We can demote an orbifold to a manifold with boundary,
and we will still have isospectrality if we impose Neumann boundary
conditions.
Dirichlet boundary conditions also work, provided that
the manifold with boundary is orientable.
This will be the case when the diagrams are treelike,
or more generally, when all cycles have even length.
If there are cycles of odd length,
when we put Dirichlet boundary conditions we must
also use twisted functions, i.e. sections of a non-trivial
bundle, which change sign when you travel around an orientation-reversing
path.

\clearpage
\newcommand{\hfig}[2]{\putfigwithsize{width=12cm}{hq#1#2}{hq#1#2}{Pair #1(#2)}}
\newcommand{\hfigsqueeze}[2]{\putfigwithsize{width=10cm}{hq#1#2}{hq#1#2}{Pair #1(#2)}}
\hfig{7}{1}
\hfig{7}{2}
\hfig{7}{3}
\clearpage
\hfig{13a}{1}
\hfig{13a}{2}
\hfig{13a}{3}
\hfig{13a}{4}
\hfig{13a}{5}
\clearpage
\hfig{13b}{1}
\hfig{13b}{2}
\hfig{13b}{3}
\hfigsqueeze{13b}{4}
\clearpage
\hfig{15}{1}
\hfig{15}{2}
\hfig{15}{3}
\hfig{15}{4}
\clearpage
\hfigsqueeze{21}{1}
\hfig{21}{2}
\hfig{21}{3}
\hfig{21}{4}
\hfig{21}{5}
\hfig{21}{6}
\hfig{21}{7}
\hfig{21}{8}
\clearpage
\hfig{11f}{1}
\hfig{11f}{2}
\hfig{11f}{3}
\hfig{11f}{4}
\clearpage
\hfig{11g}{1}
\hfig{11g}{2}
\hfig{11g}{3}
\hfig{11g}{4}
\hfig{11g}{5}
\hfig{11g}{6}
\clearpage
\hfig{11h}{1}
\hfig{11h}{2}
\hfigsqueeze{11h}{3}
\hfig{11h}{4}
\clearpage
\hfig{11i}{1}
\hfig{11i}{2}
\hfig{11i}{3}
\hfig{11i}{4}
\hfig{11i}{5}

\clearpage
\appendix
\section{The master at work} \label{appendix}
Here,
in the hand of the master,
are quilts,
along with associated diagrams for isospectral pairs
(including peacocks rampant and couchant
in the hand of the pupil).

\vspace{2.0cm}

\includegraphics[width=\textwidth]{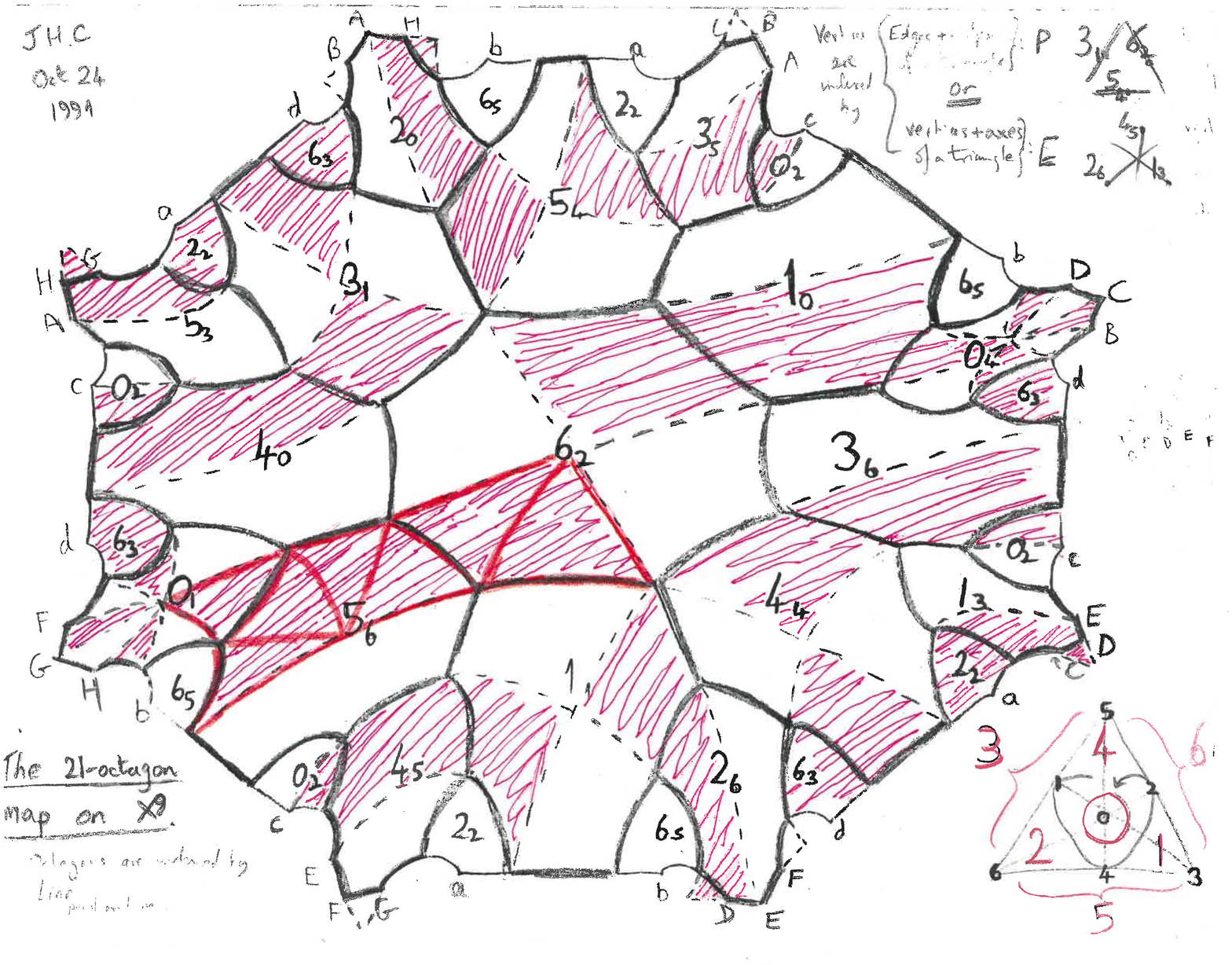}
\centerline{The 21-octagon map on $\cross^9$.
Octagons are indexed by $\mbox{\small line}_{\mbox{\tiny point on line}}$.
}

\newcommand{\figpage}[1]{
\includegraphics[height=.95\textheight]{scan/#1.pdf}
}

\figpage{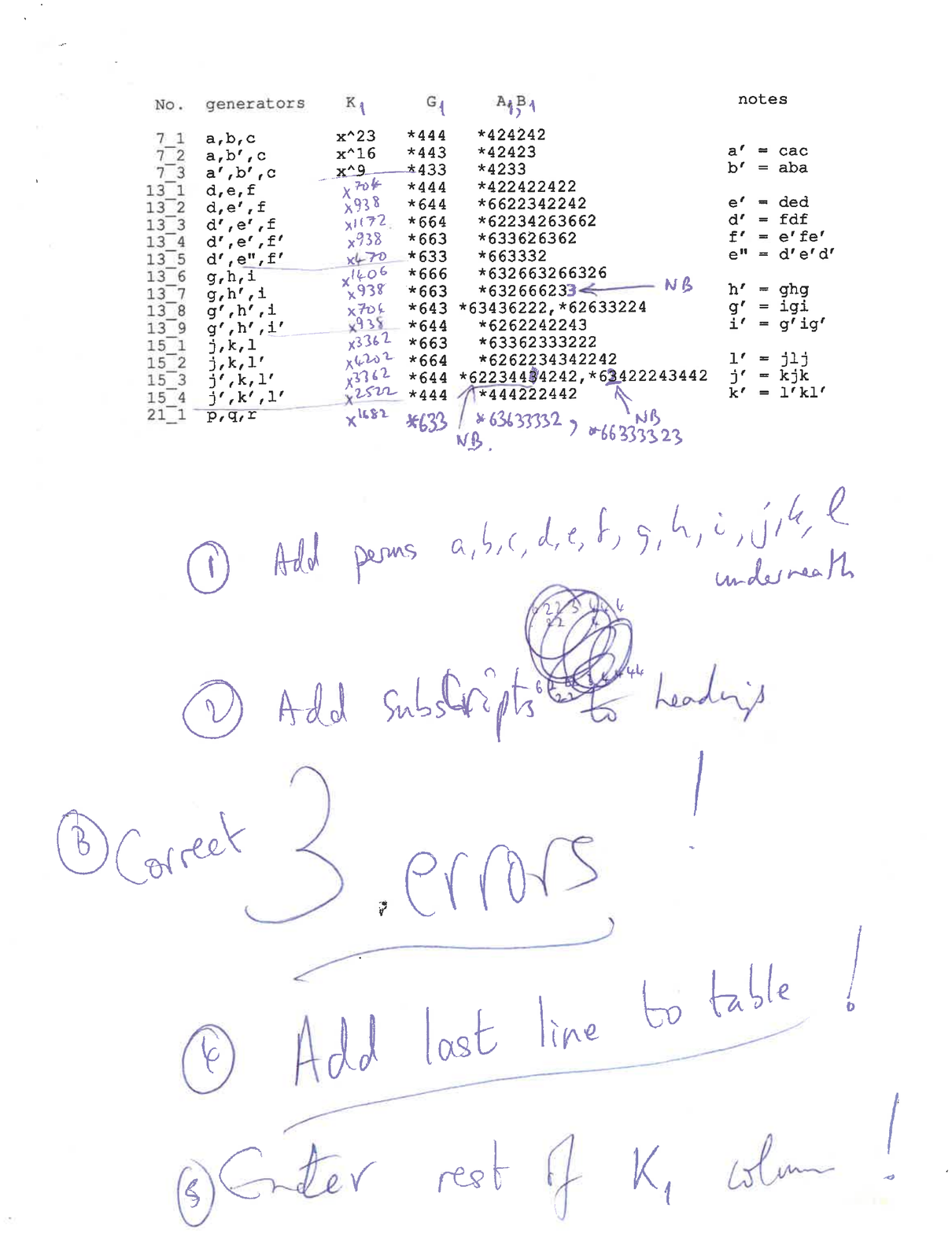}

\figpage{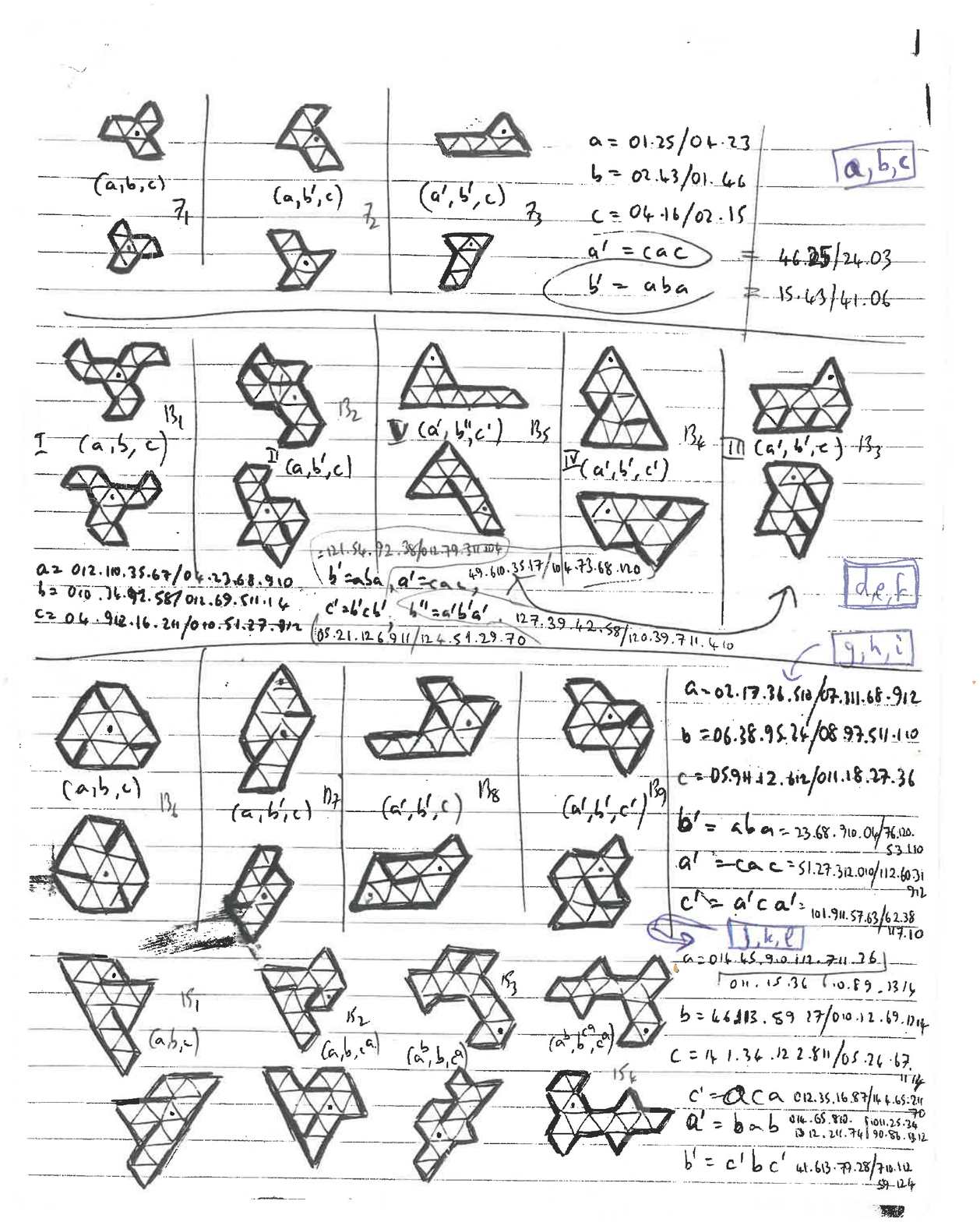}

\figpage{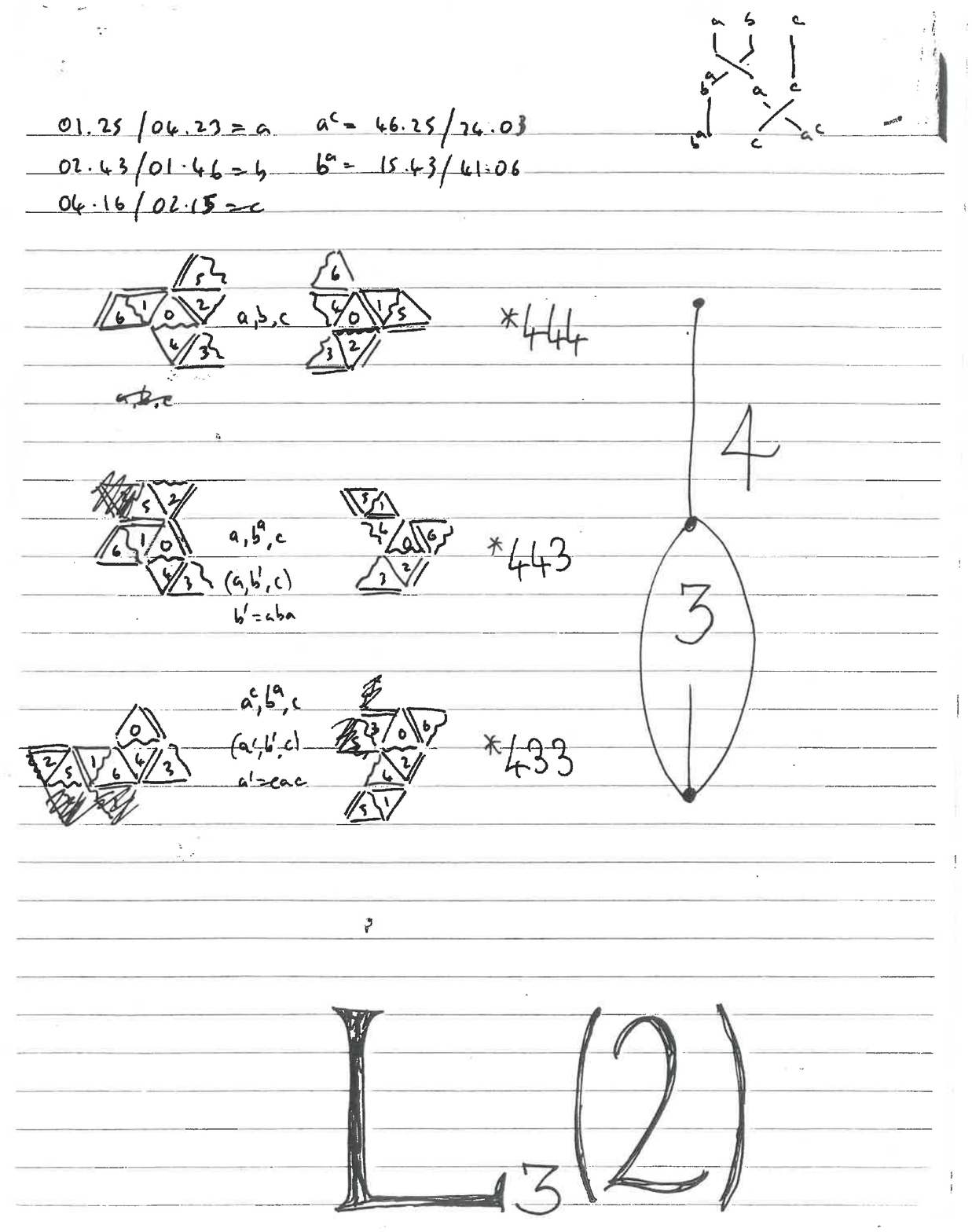}

\figpage{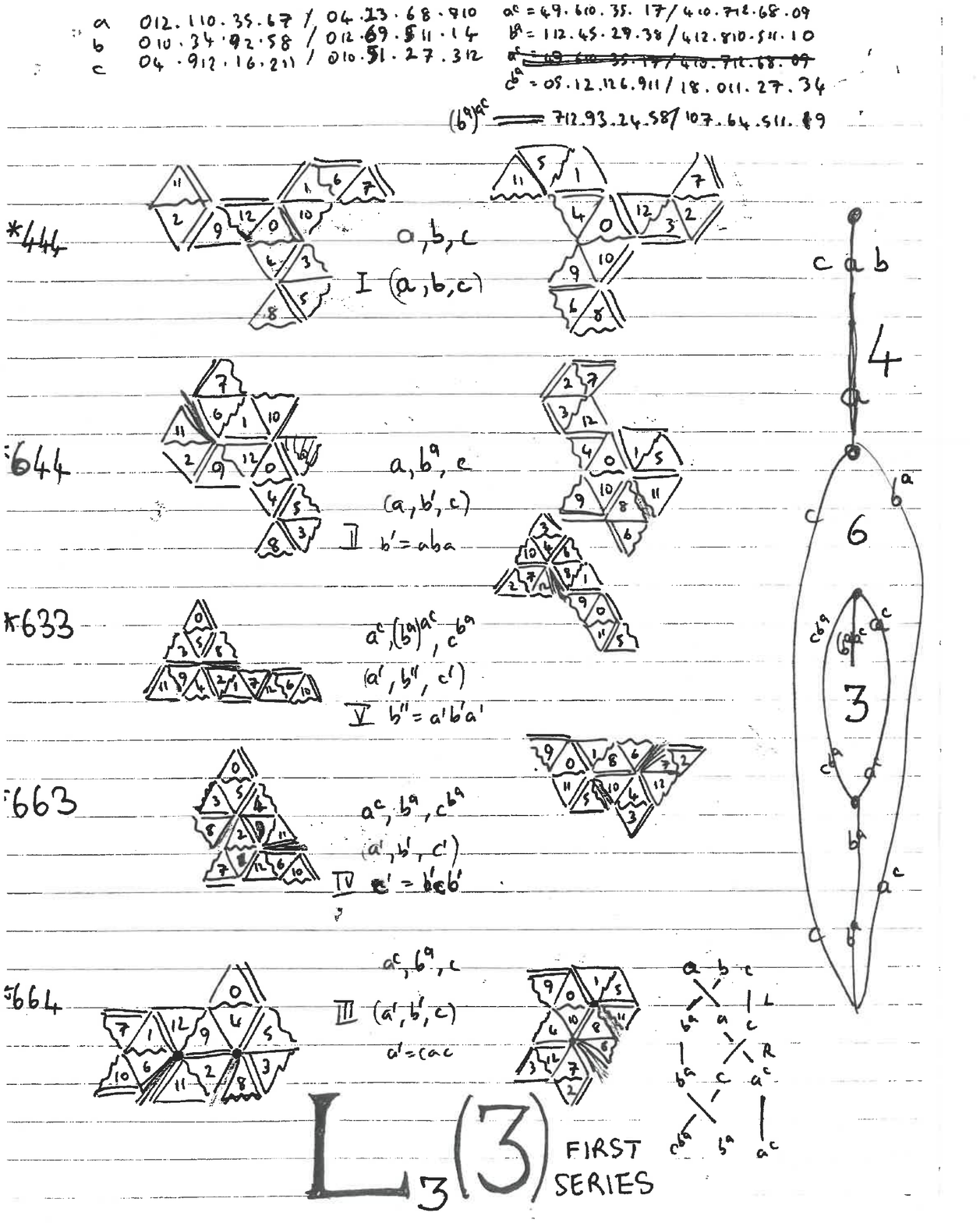}

\figpage{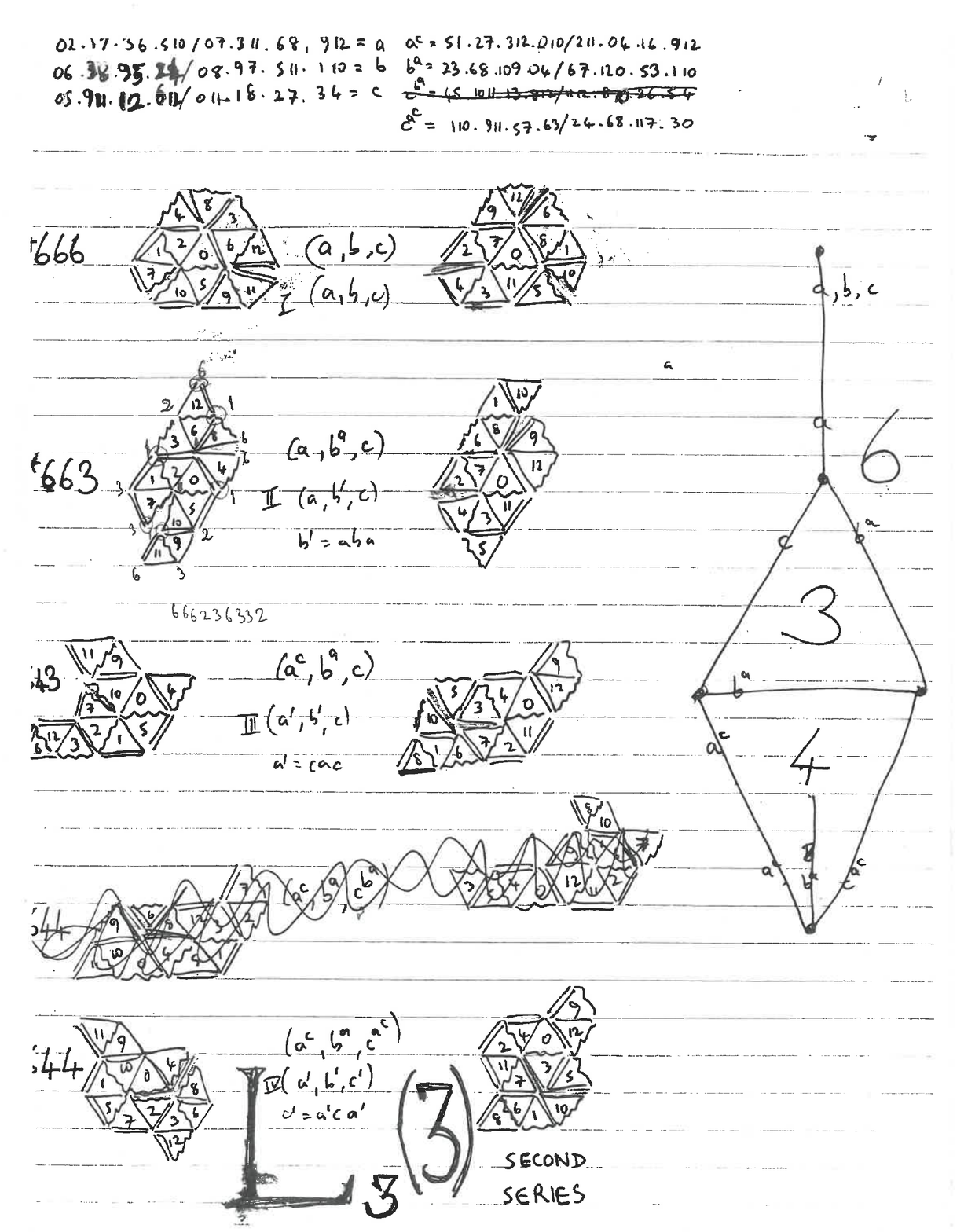}

\figpage{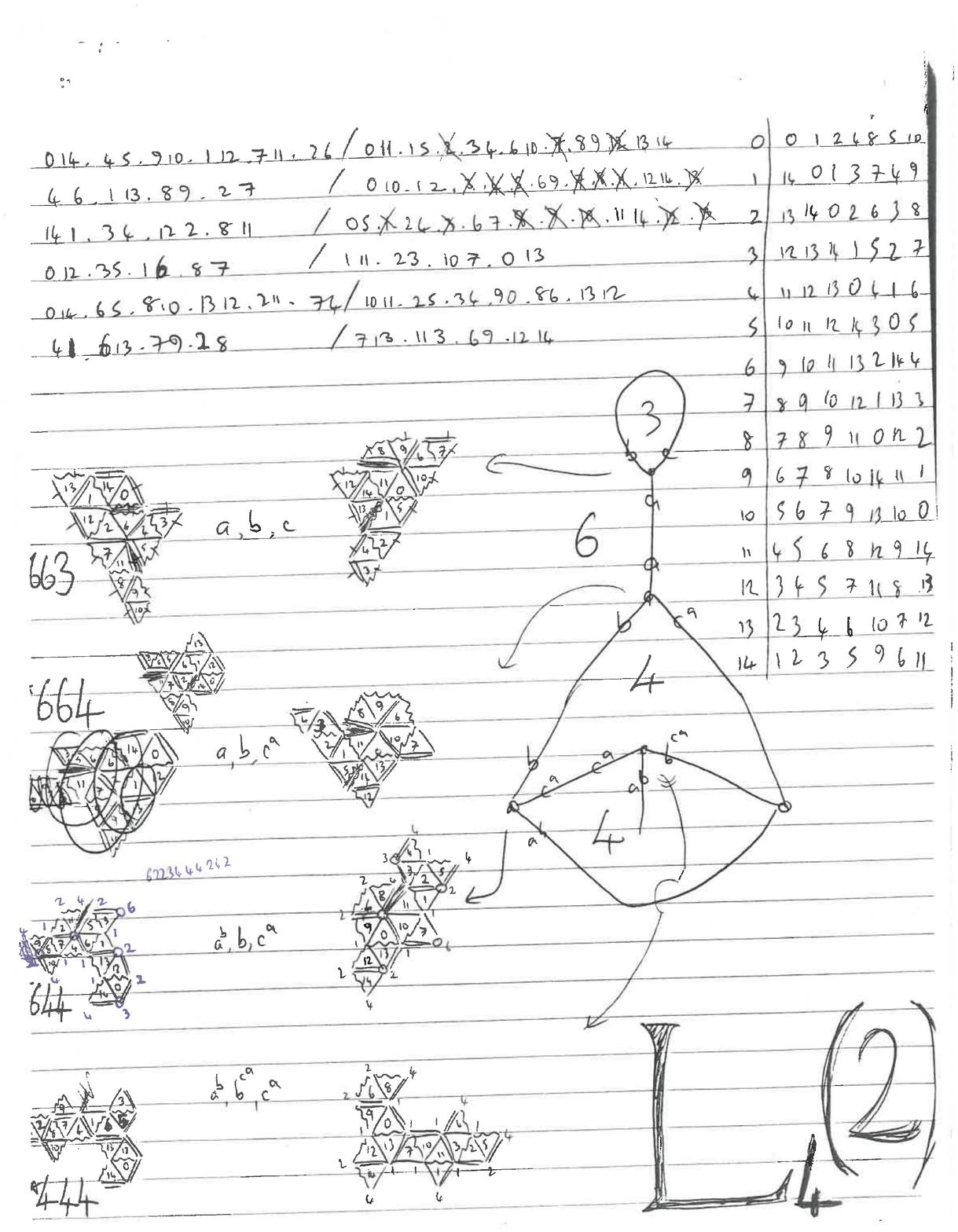}

\figpage{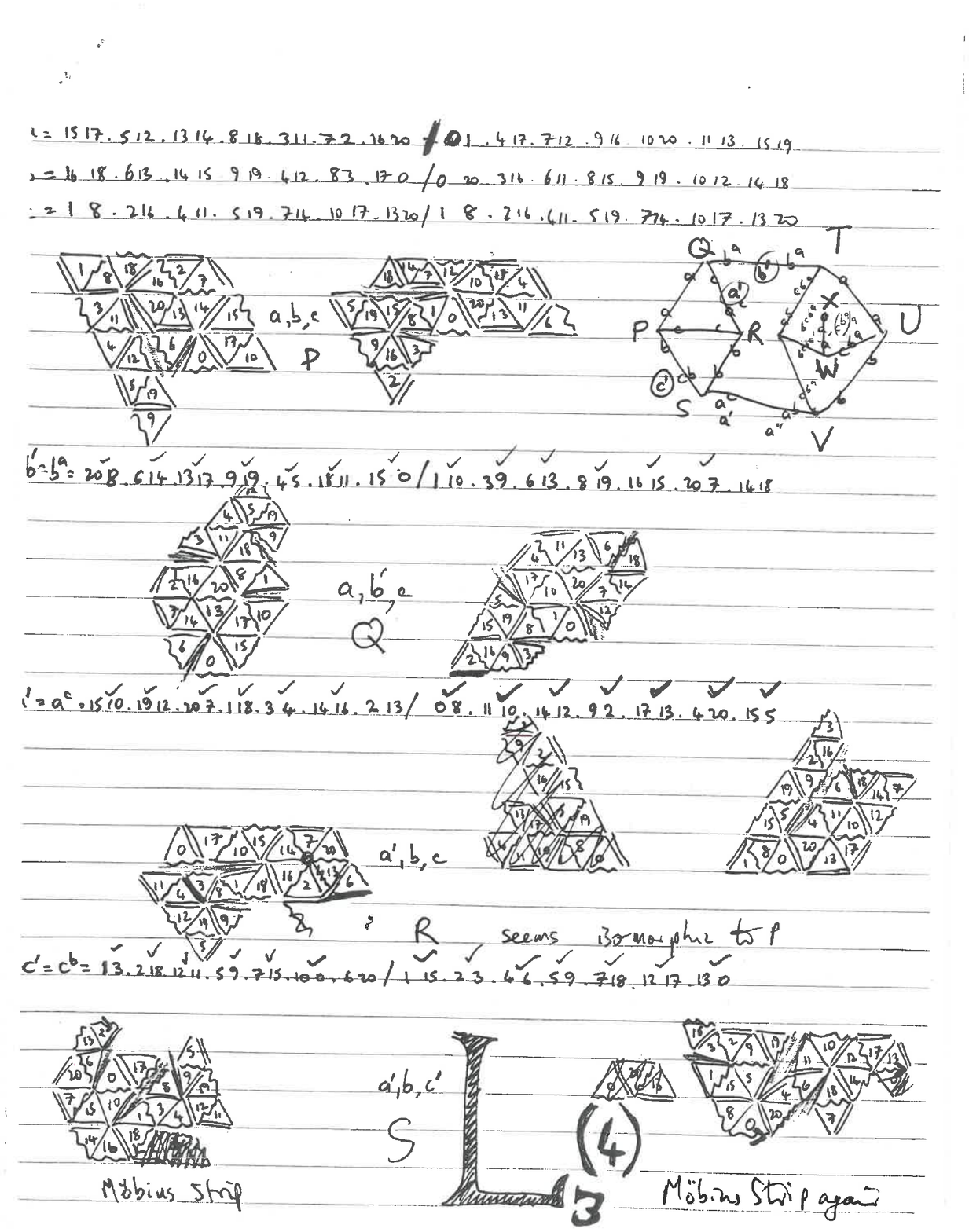}

\figpage{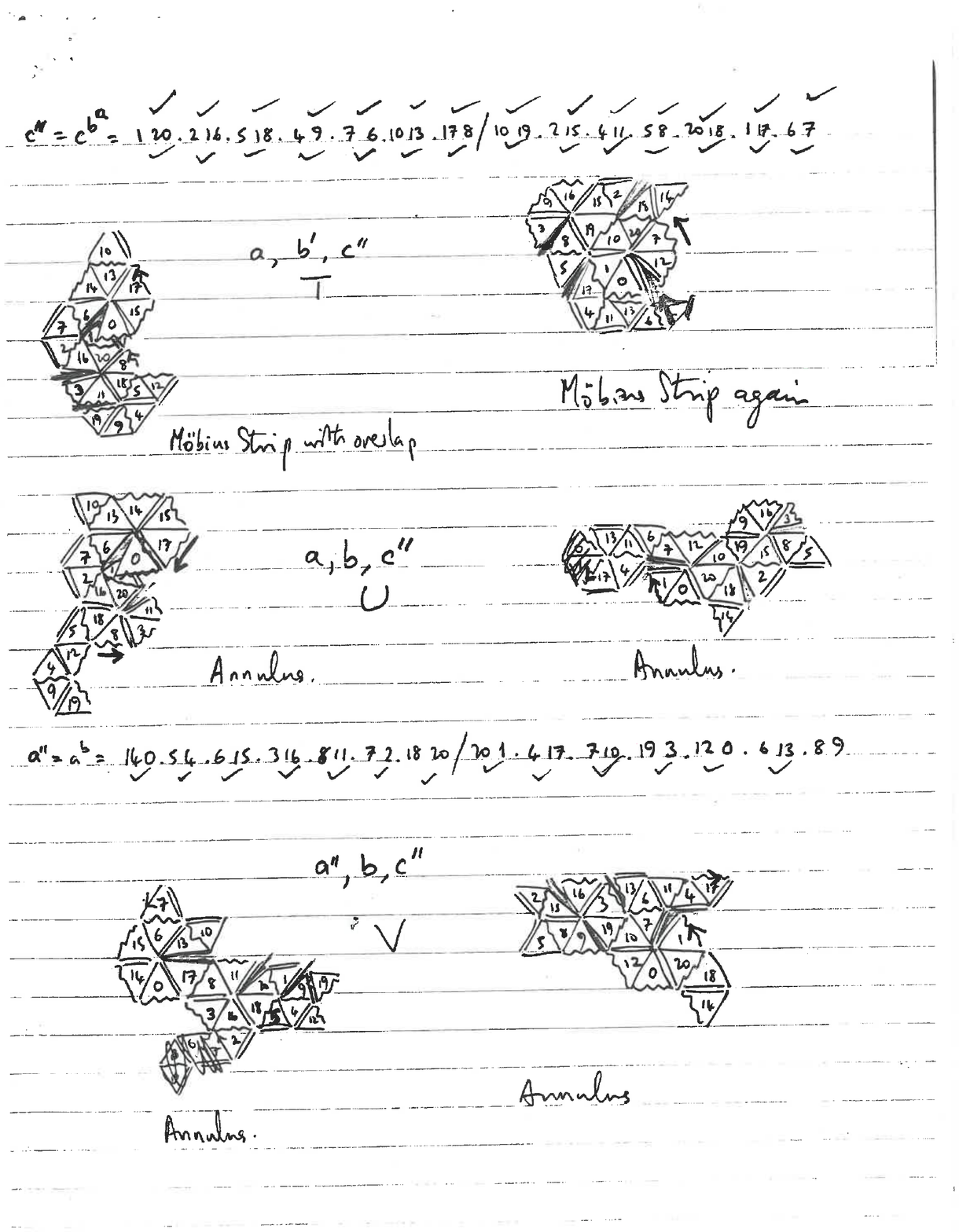}

\figpage{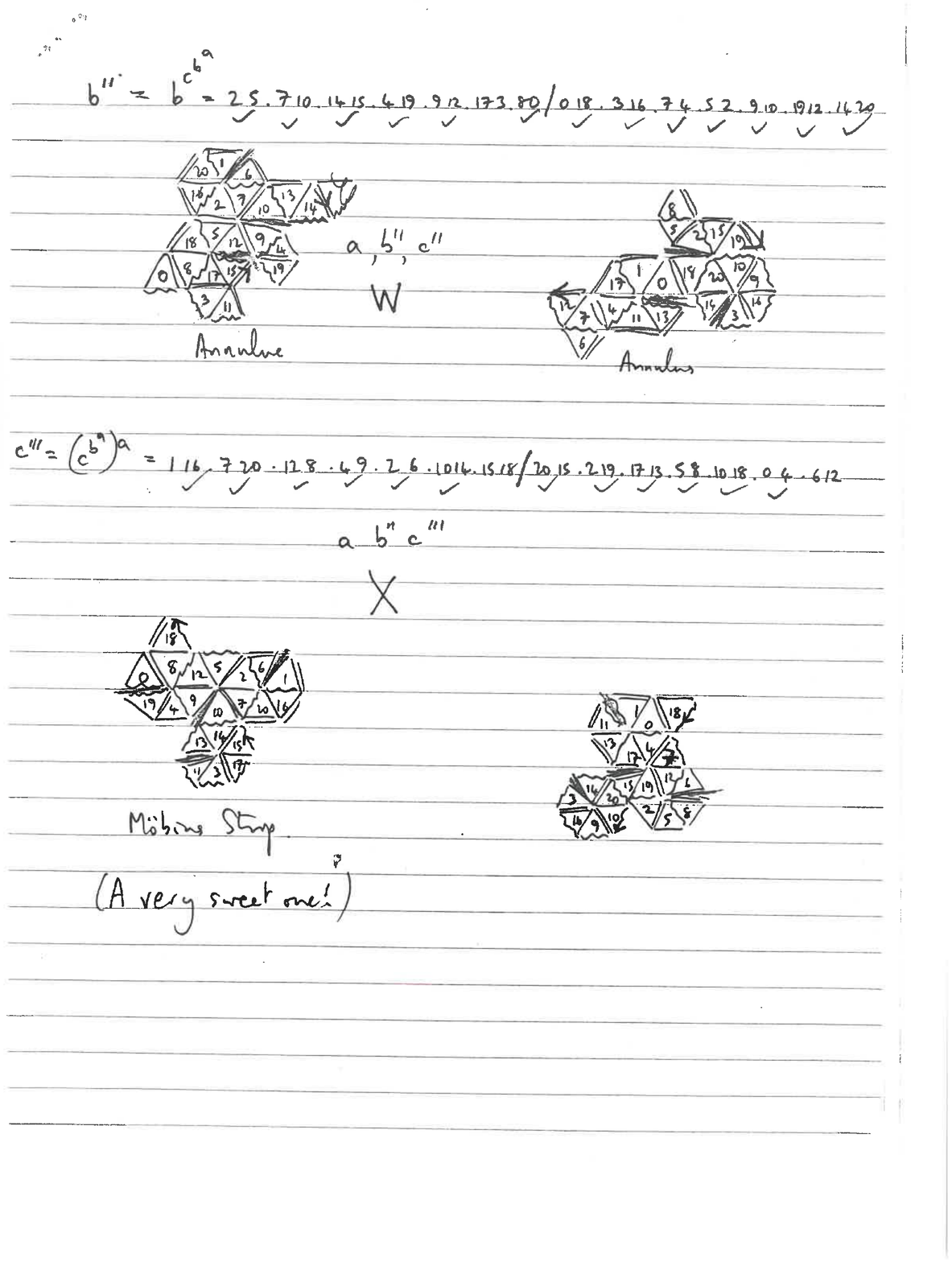}




\figpage{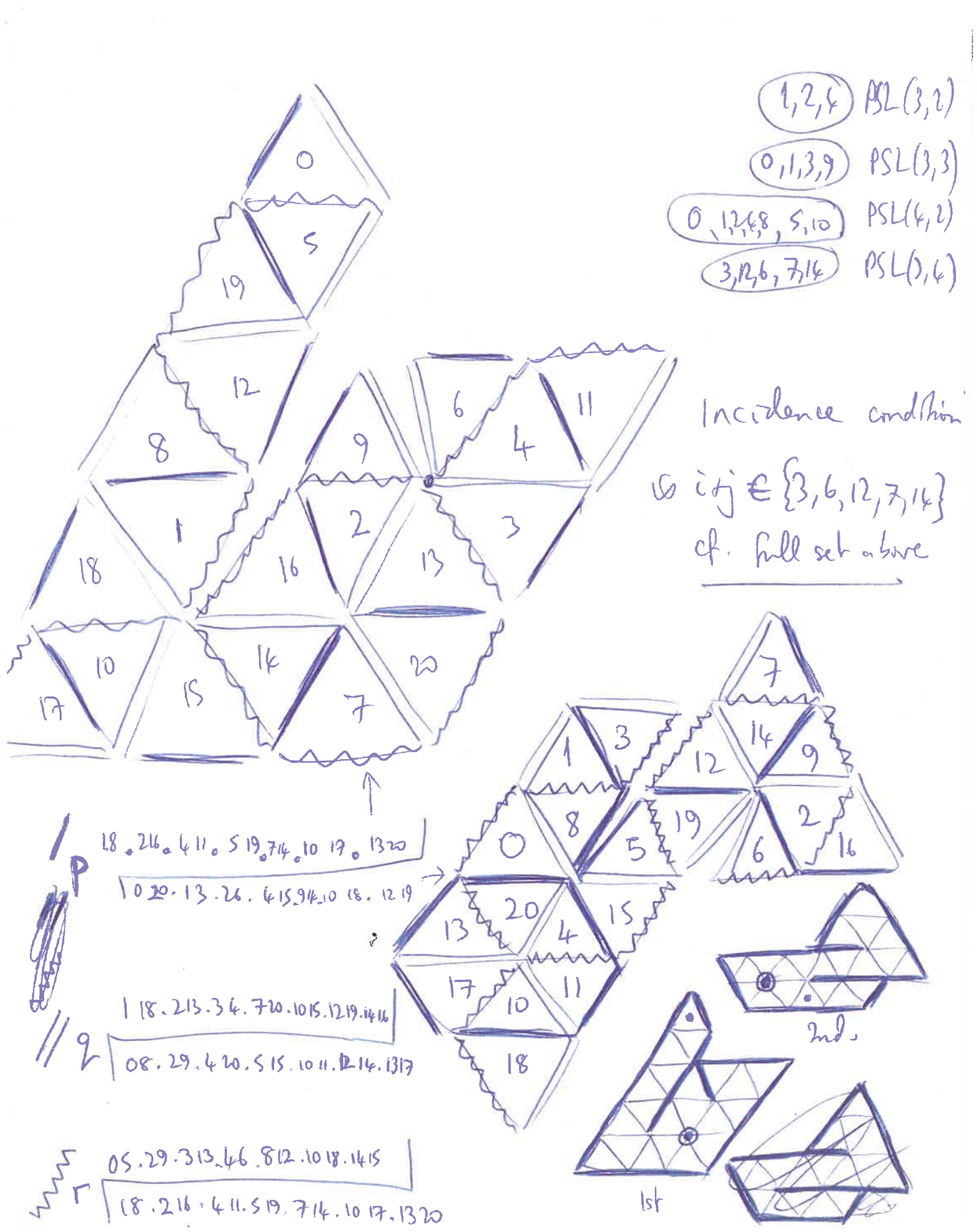}

\figpage{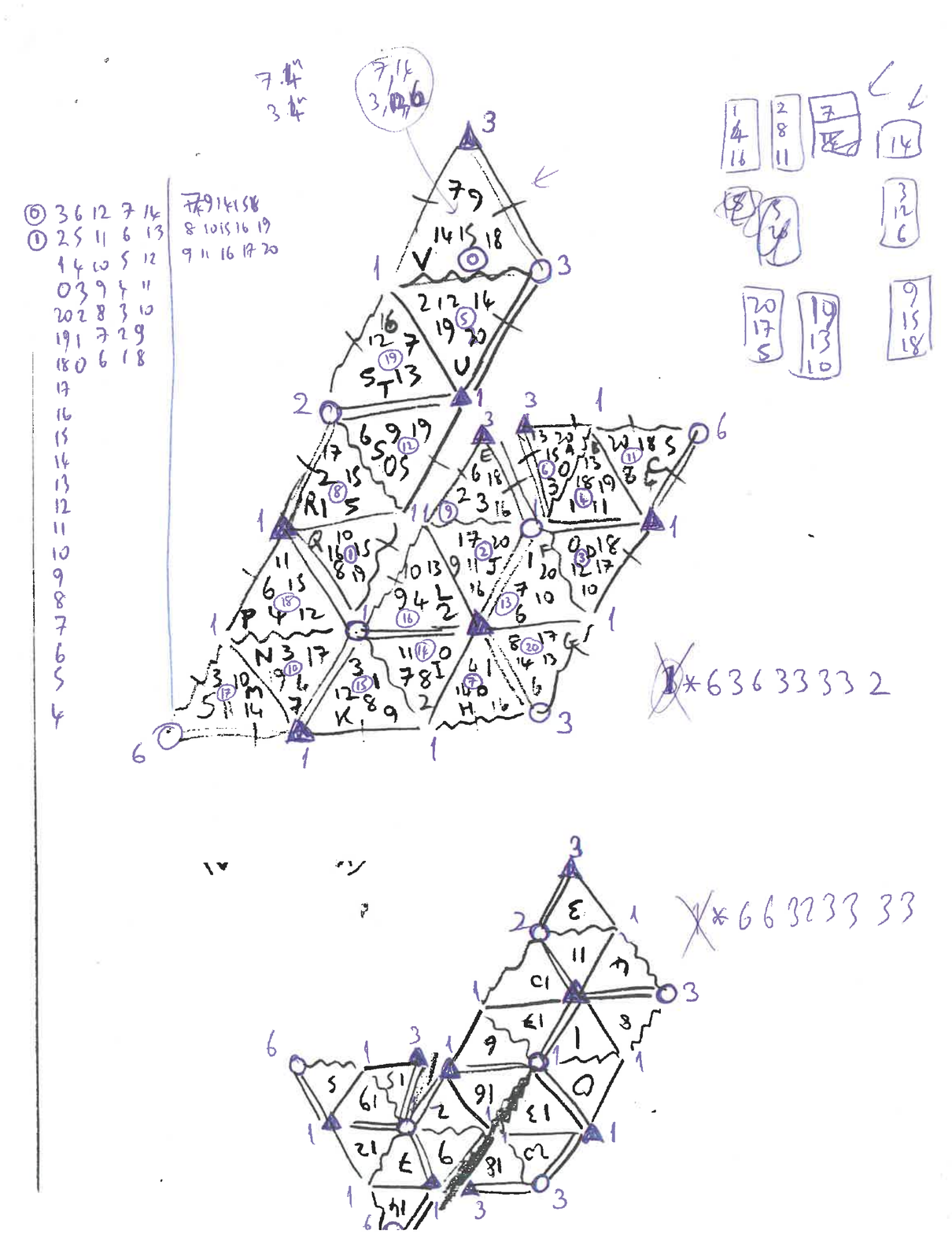}


\figpage{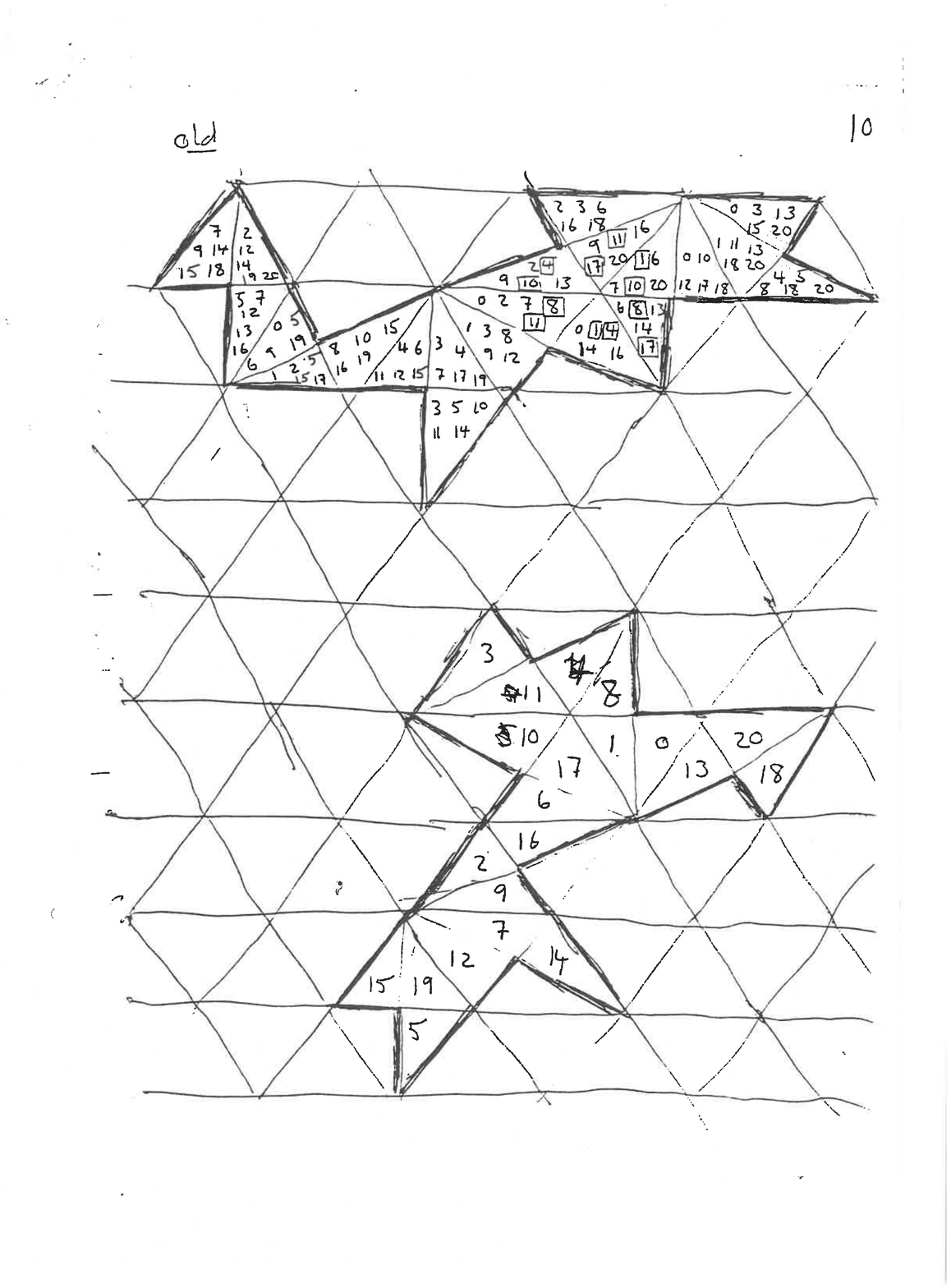}

\figpage{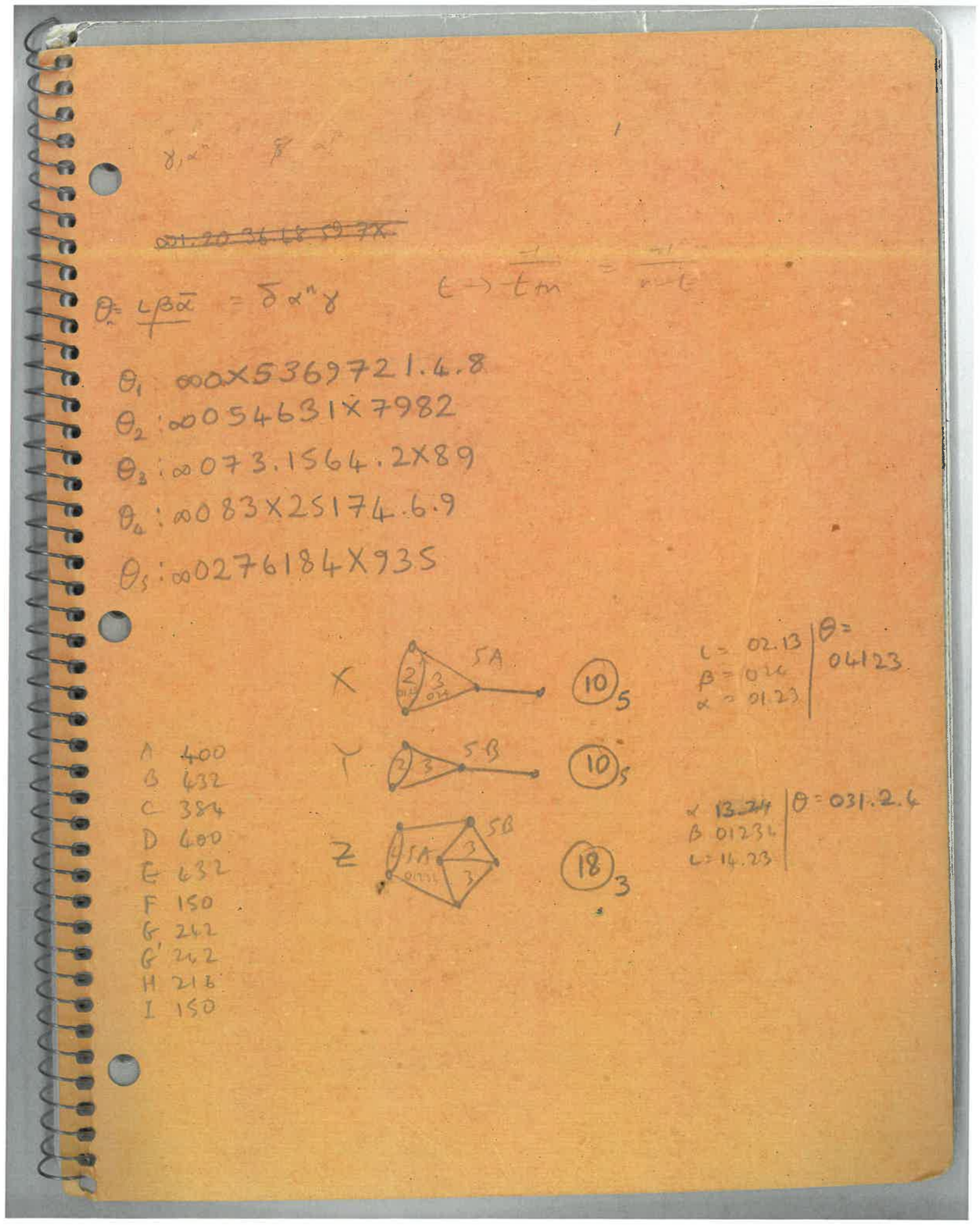}

\figpage{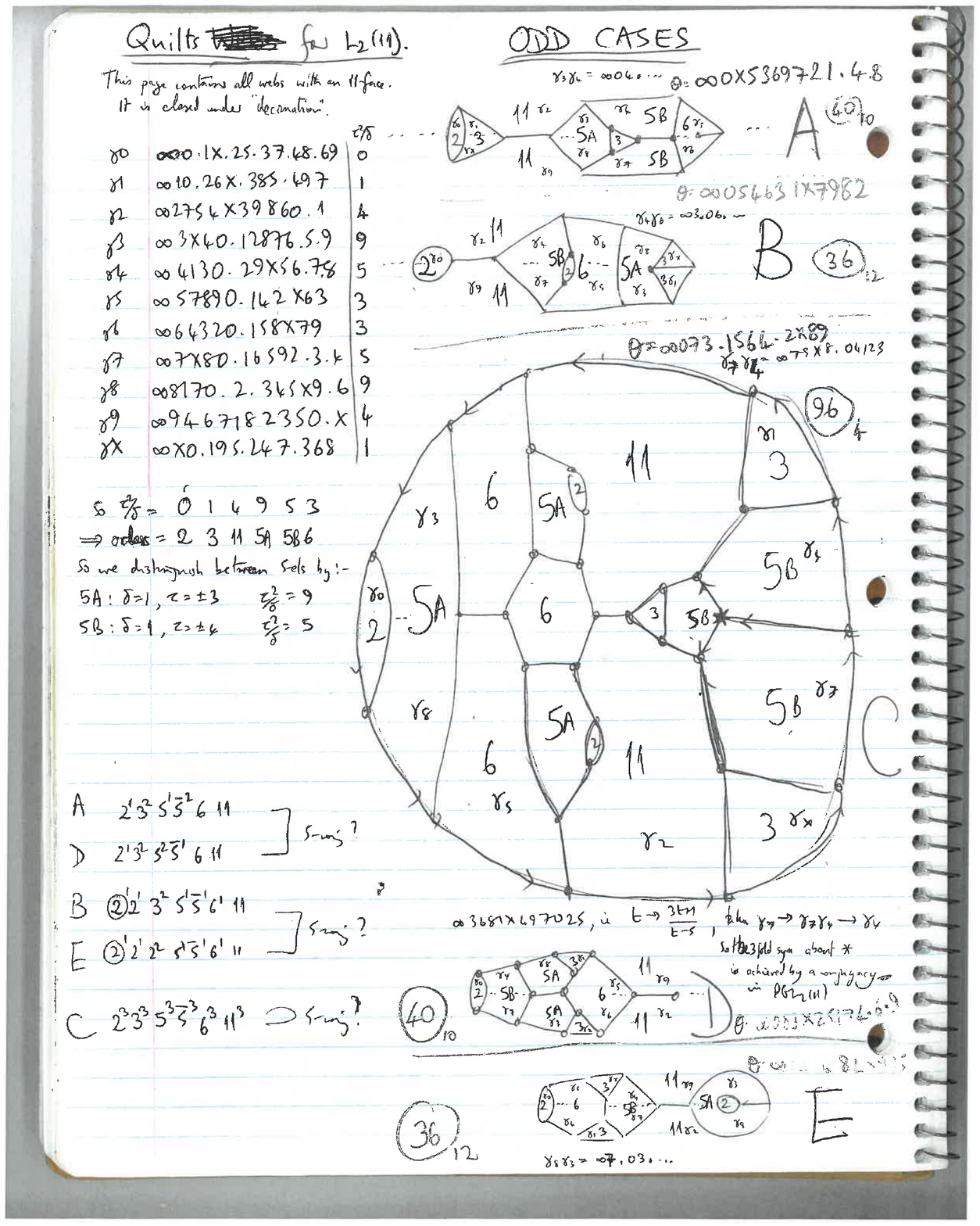}

\figpage{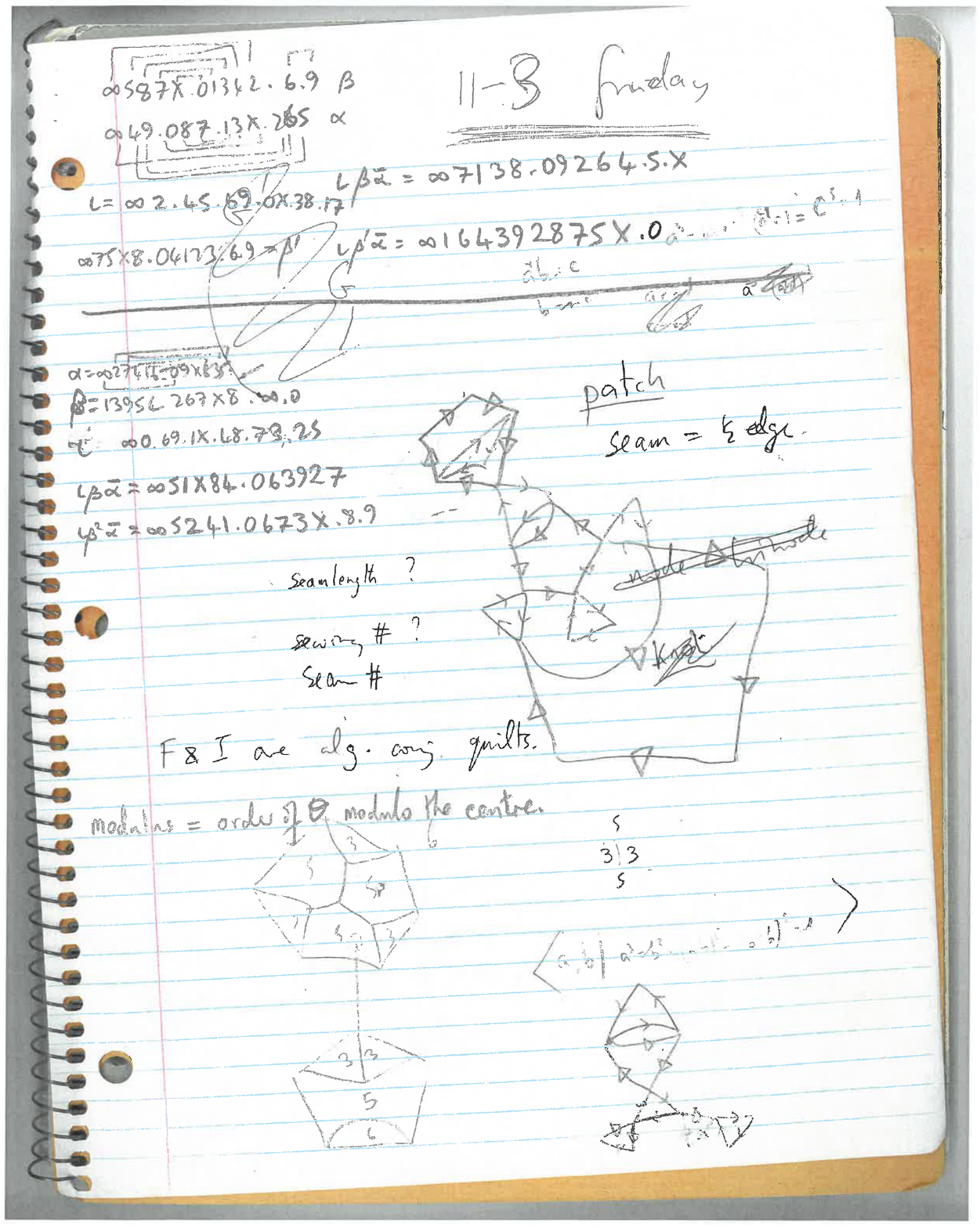}

\figpage{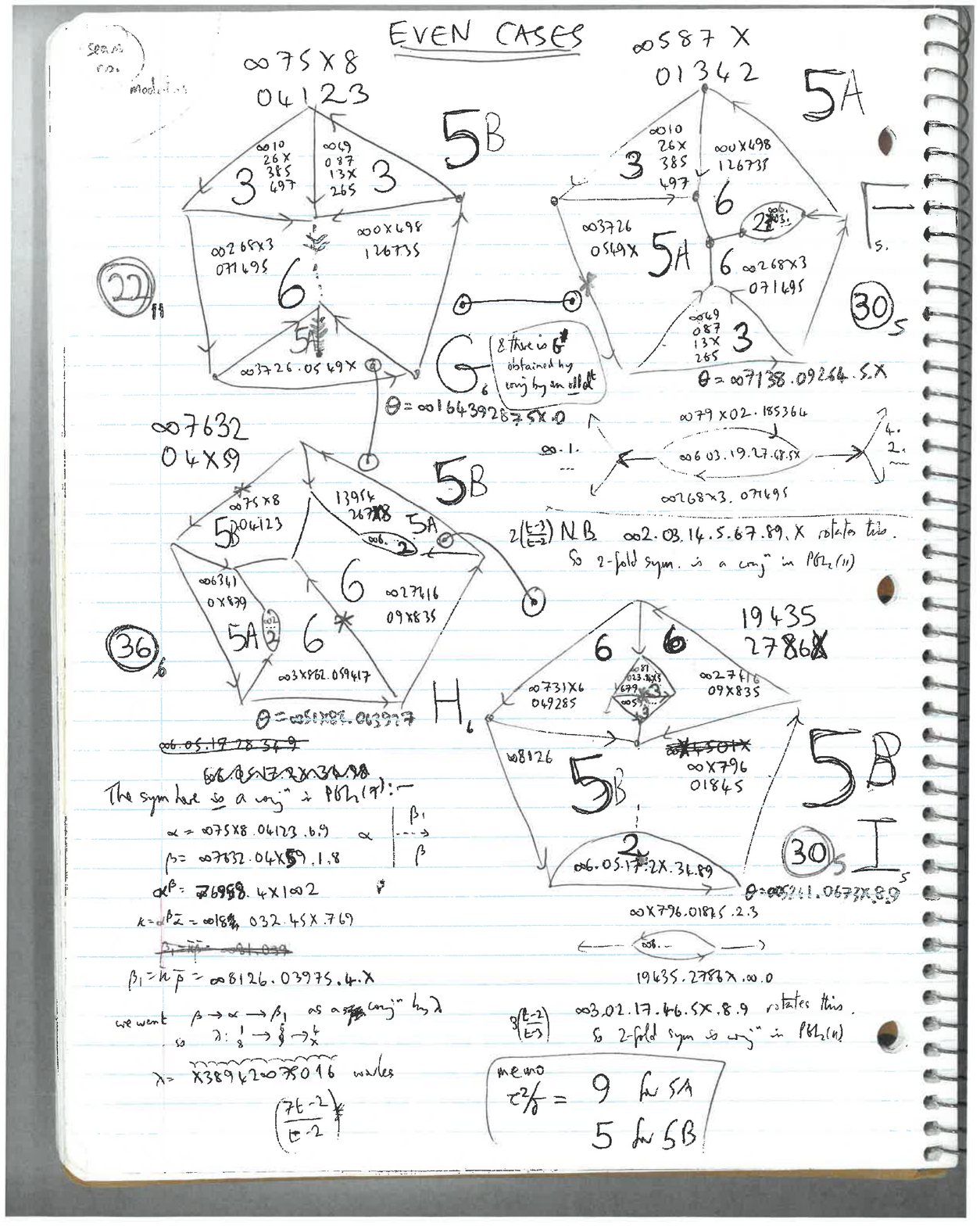}

\figpage{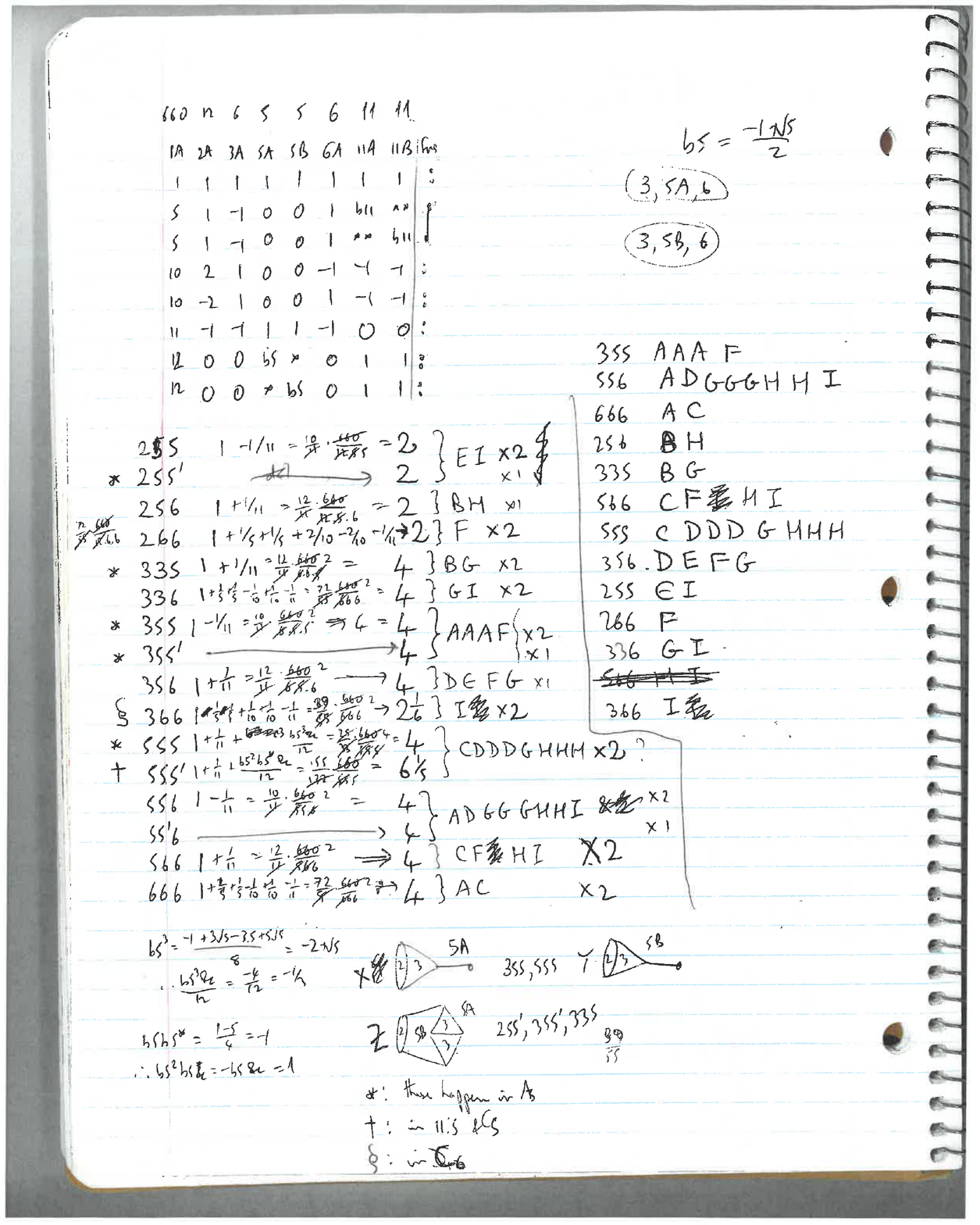}

\figpage{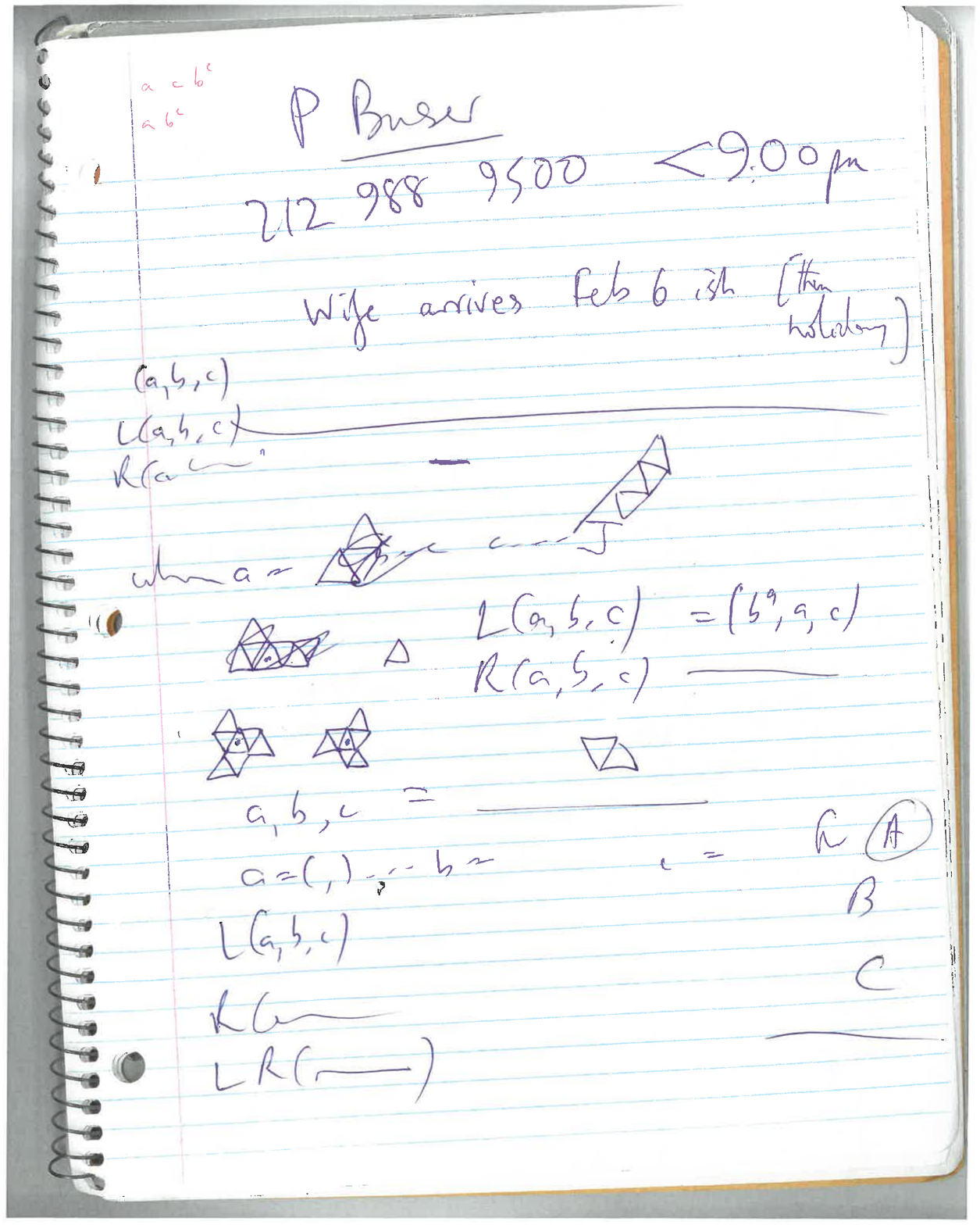}

\figpage{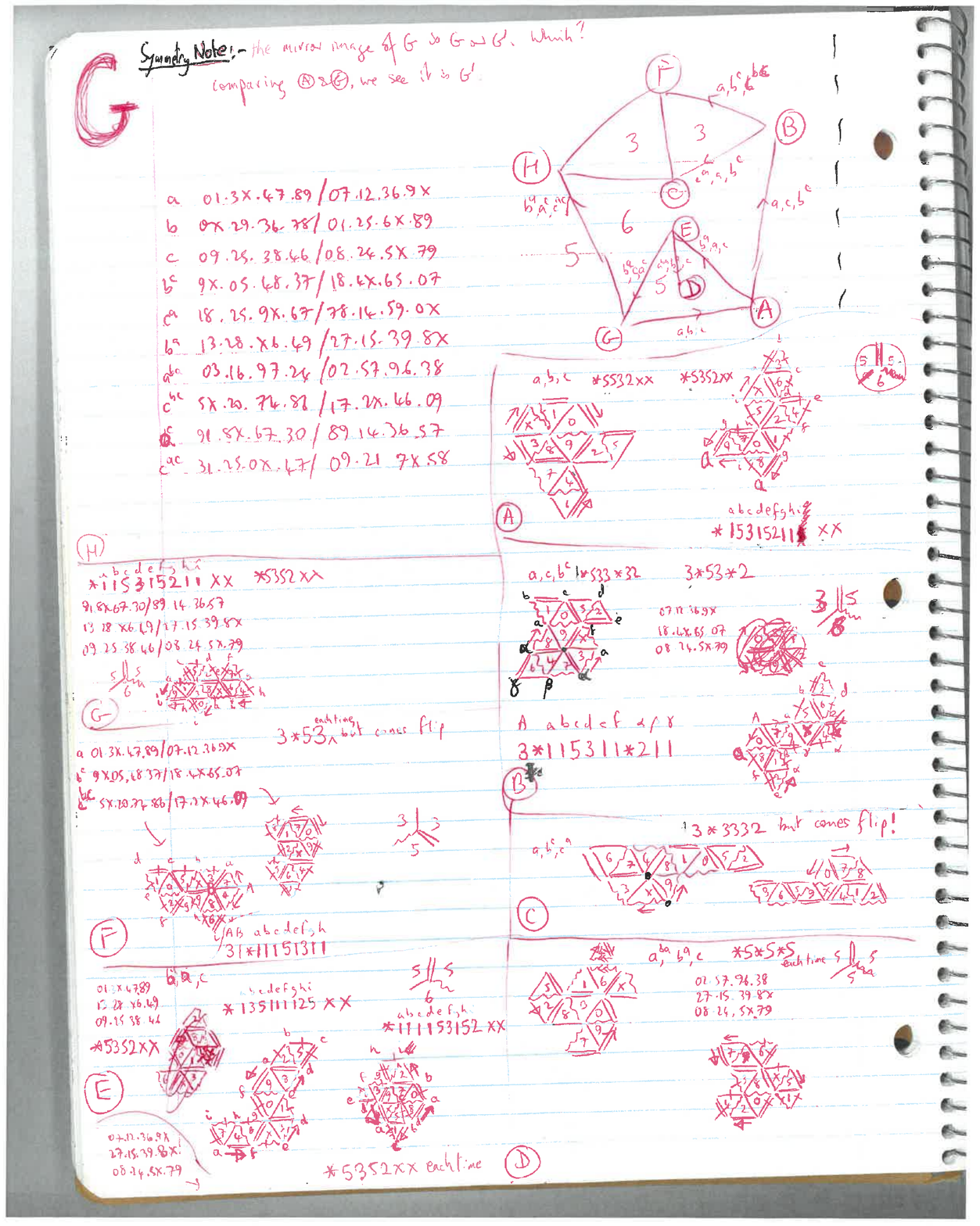}

\figpage{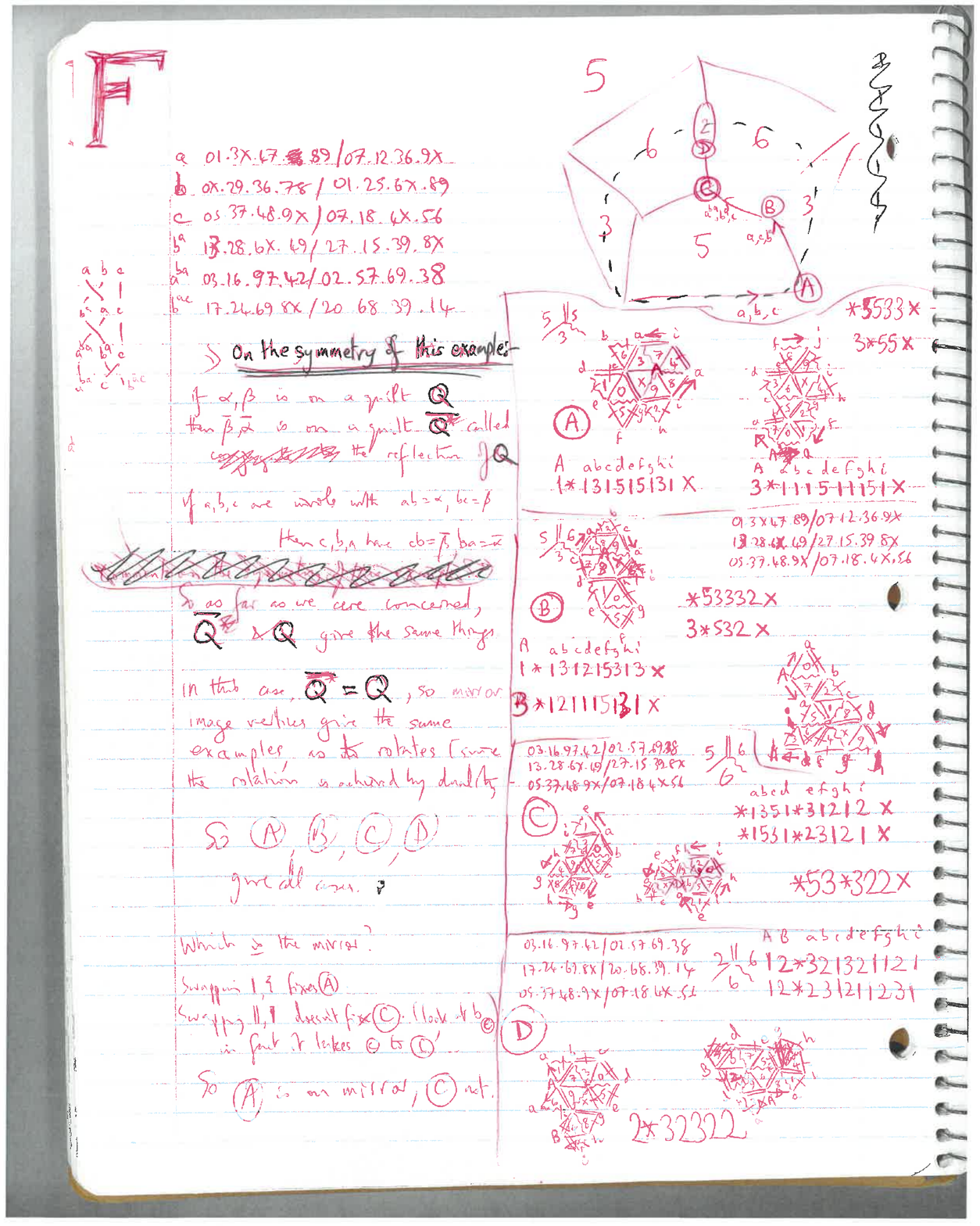}

\figpage{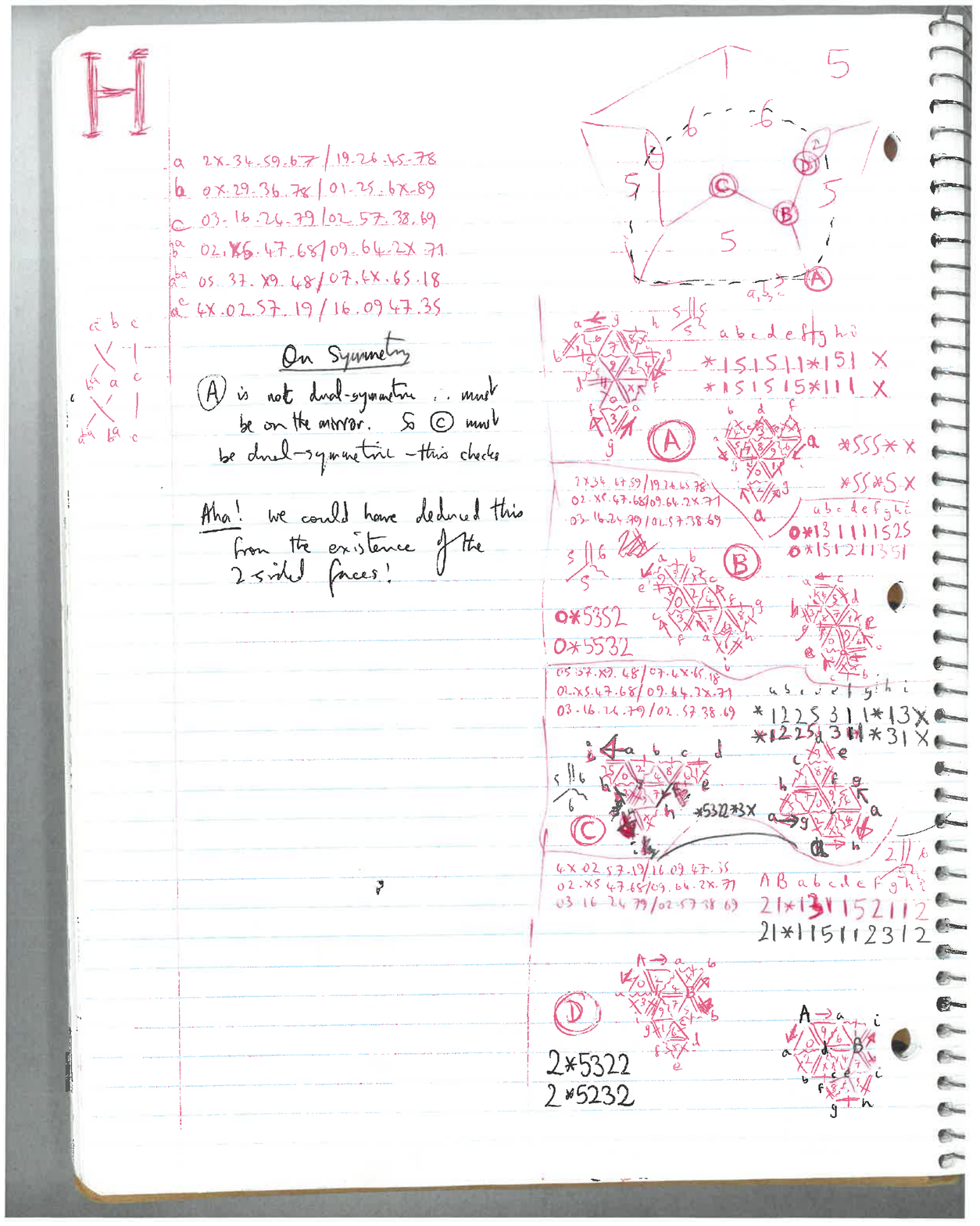}

\figpage{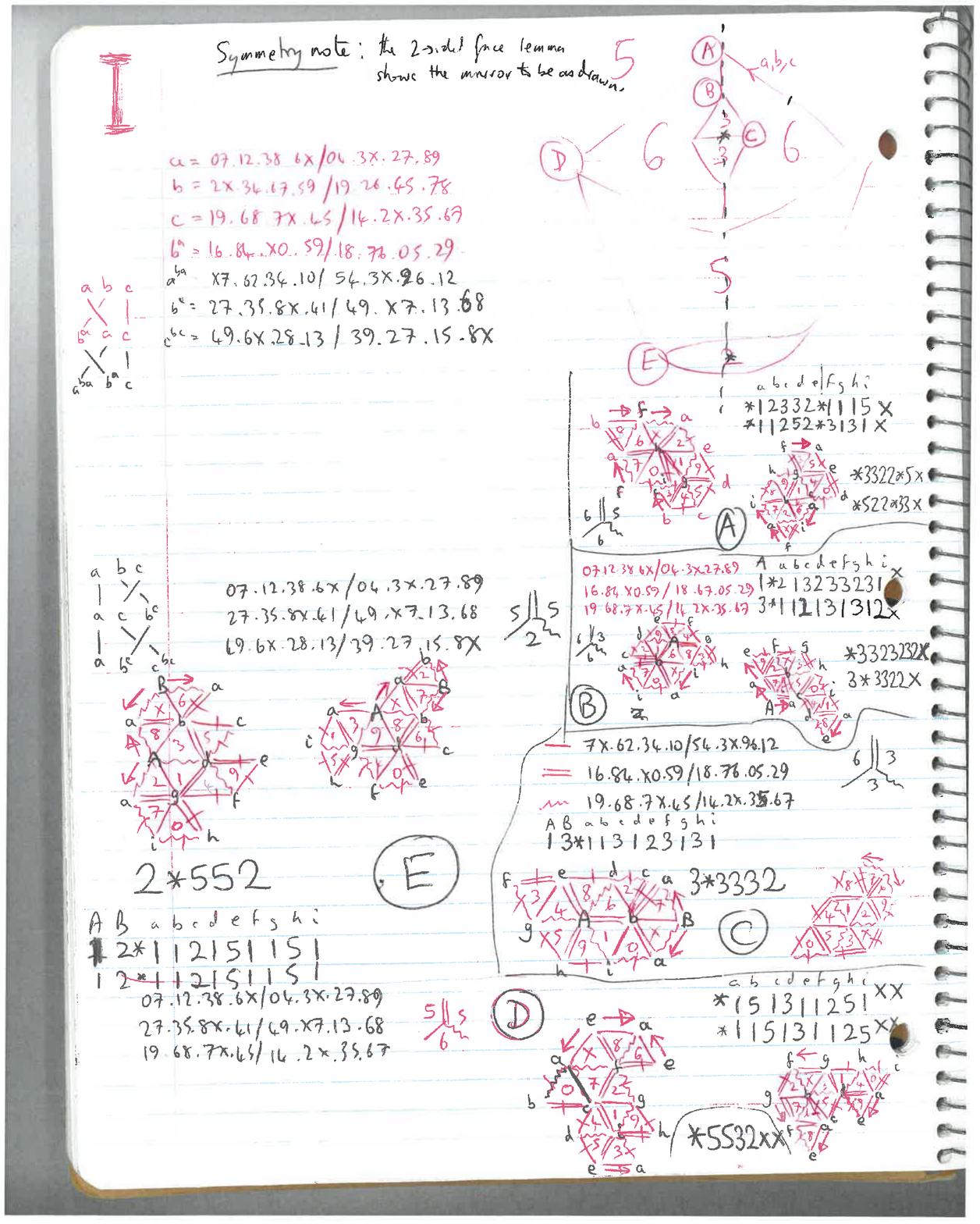}

\figpage{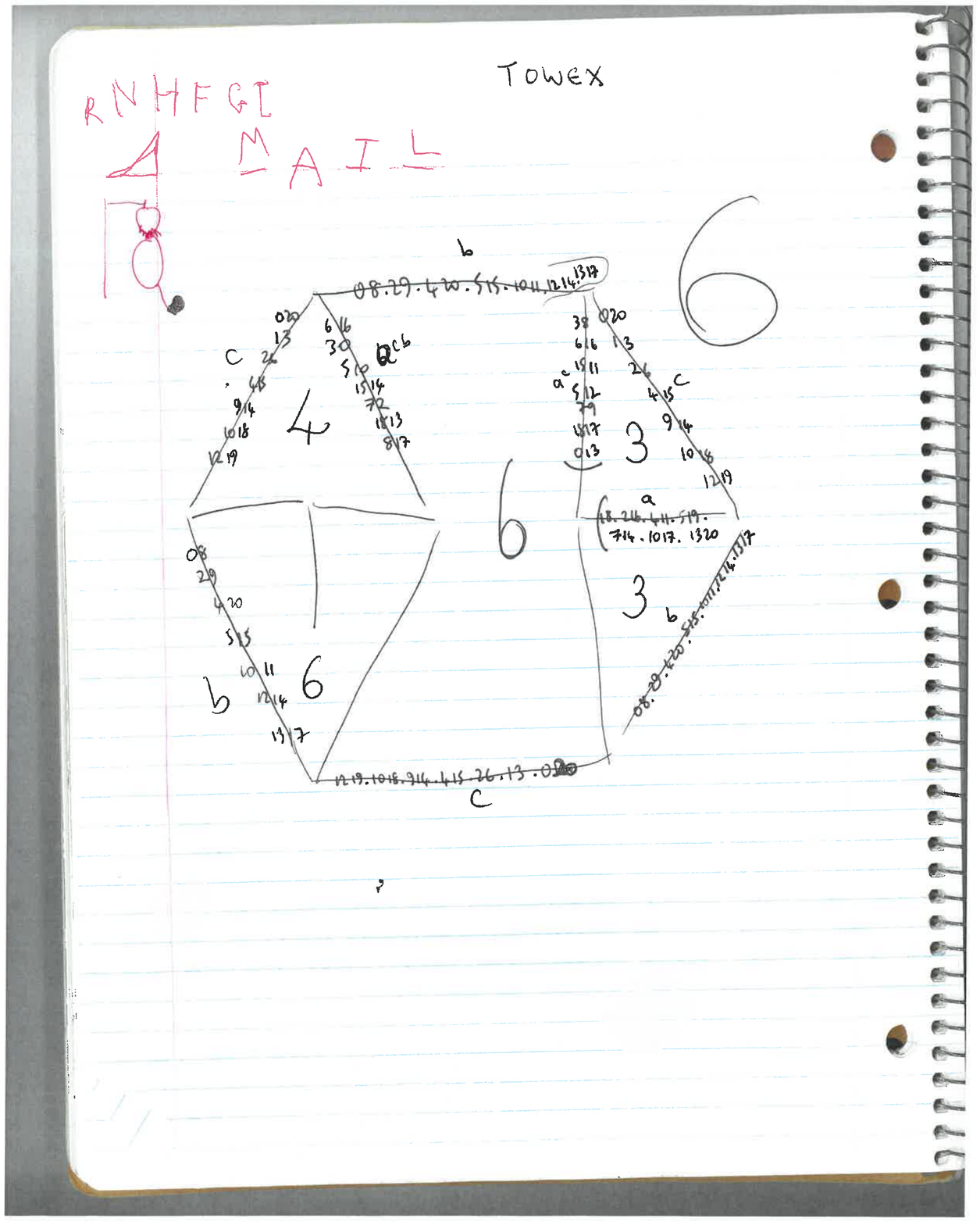}

\clearpage
\bibliography{quilt}
\bibliographystyle{hplain}
\end{document}